%% file: main.tex
\documentclass[11pt]{article}
\usepackage{amssymb,amsmath,amsthm}
\usepackage[colorlinks]{hyperref}
\usepackage[all]{xy}

\title{New derived symmetries \\ of some hyperk\"ahler varieties}
\author{Nicolas Addington}
\date{}

\newcommand \A {\mathcal A}
\newcommand \B {\mathcal B}
\newcommand \C {\mathbb C}
\newcommand \E {\mathcal E}
\newcommand \F {\mathcal F}
\newcommand \G {\mathcal G}
\renewcommand \H {\mathcal H}
\newcommand \I {\mathcal I}

\newcommand \M {\mathcal M}
\renewcommand \O {\mathcal O}
\renewcommand \P {\mathbb P}
\newcommand \Q {\mathbb Q}
\newcommand \U {\mathcal U}
\newcommand \Z {\mathbb Z}

\DeclareMathOperator \id {id}
\DeclareMathOperator \Aut {Aut}
\DeclareMathOperator \Ext {Ext}
\DeclareMathOperator \RHom {RHom}
\DeclareMathOperator \cone {cone}
\DeclareMathOperator \Spec {Spec}
\DeclareMathOperator \Hom {Hom}
\DeclareMathOperator \im {im}
\DeclareMathOperator \RGamma {R\Gamma}
\DeclareMathOperator \Bl {Bl}
\DeclareMathOperator \Proj {Proj}
\DeclareMathOperator \Sym {Sym}
\DeclareMathOperator \supp {supp}
\DeclareMathOperator \Gr {Gr}
\DeclareMathOperator \Tor {Tor}
\DeclareMathOperator \rank {rank}
\DeclareMathOperator \ev {eval}
\newcommand \HH {H\!H}

\newtheorem*{conj*}{Conjecture}
\newtheorem{thm}{Theorem}
\newtheorem{prop}{Proposition}[section]
\newtheorem{lem}[prop]{Lemma}
\theoremstyle{definition}
\newtheorem{defn}[prop]{Definition}

\numberwithin{equation}{section}

\begin{document}
\maketitle
\input intro
\input spherical

\input hilb_calc
\input P-functors

\input cubic_calc

\input base_change

\pagebreak
\newcommand \httpurl [1] {\href{http://#1}{\nolinkurl{#1}}}
\bibliographystyle{plain}
\bibliography{main}

\end{document}

%% file: intro.tex

\begin{abstract}
We construct a new autoequivalence of the derived category of the Hilbert scheme of $n$ points on a K3 surface, and of the variety of lines on a smooth cubic 4-fold.  For Hilb$^2$ and the variety of lines, we use the theory of spherical functors; to deal with Hilb$^n$ for $n>2$ we develop a theory of $\P$-functors.  We conjecture that the same construction yields an autoequivalence for any moduli space of sheaves on a K3 surface.

In an appendix we give a cohomology and base change criterion which is well-known to experts, but not well-documented.
\end{abstract}

\section*{Introduction}

This paper grows out of the following observation:  Let $S$ be a complex K3 surface, let $S^{[2]}$ be the Hilbert scheme of pairs of points on $S$, thought of as a moduli space of ideal sheaves, let $F\colon D^b(S) \to D^b(S^{[2]})$ be the functor induced by the universal sheaf $\U$ on $S \times S^{[2]}$, and let $R$ be the right adjoint of $F$.  Then the composition $RF$ is isomorphic to $\id_S \oplus [-2]$, so $F$ is a ``spherical functor'' in the sense of Rouquier \cite{rouquier} and Anno \cite{anno} and hence determines an auto\-equivalence $T$ of $D^b(S^{[2]})$.  Briefly, $T = \cone(FR \to \id)$.  Spherical functors generalize Seidel and Thomas's spherical objects \cite{st} and unify various family versions of them \cite{horja,toda}.  In \S\ref{spherical} we give a simplified definition of spherical functors, review the known examples, and give an alternate proof that they yield autoequivalences in preparation for our work on $\P$-functors in \S\ref{P-functors}.

The Fourier--Mukai kernel inducing $T$ is a shift of the sheaf
\[ \E xt^1_{\pi_{13}}(\pi_{12}^* \U, \pi_{23}^* \U) \]
on $S^{[2]} \times S^{[2]}$, where $\pi_{ij}$ are the projections from $S^{[2]} \times S \times S^{[2]}$.  Markman has studied this sheaf in his paper \cite{markman} on the Beauville--Bogomolov form.  It is a reflexive sheaf of rank 2, locally free away from the diagonal.  Thus $T$ sends the structure sheaf of a point to a sheaf of rank 2, so it is not in the subgroup of $\Aut(D^b(S^{[2]}))$ generated by shifts, line bundles, automorphisms of $S^{[2]}$, and $\P$-twists, all of which preserve rank (up to sign).  We will also see (\S\ref{spherical_coho}) that it does not come from any known spherical twist on $S$ via the map $\Aut(D^b(S)) \hookrightarrow \Aut(D^b(S^{[2]}))$ studied by Ploog \cite{ploog_paper}.

Next we ask what happens when we replace $S^{[2]}$ with $S^{[n]}$.  In \S\ref{hilb_calc} we show that
\[ RF = \id_S \oplus [-2] \oplus [-4] \oplus \dotsb \oplus [-2n+2]. \]
Markman and Mehrotra \cite[Thm.~2.2(1)]{mm} have given another proof of this using the equivalence $D^b(S^{[n]}) \cong D^b([S^n/\mathfrak S_n])$ of Bridgeland--King--Reid and Haiman; our proof is more geometric.  To get an autoequivalence of $D^b(S^{[n]})$ from our functor $F$, we are obliged to generalize Huybrechts and Thomas's $\P$-objects \cite{ht}.  In \S\ref{P-functors} we define $\P$-functors, give more examples of them, and show that they yield new autoequivalences.

The behavior we are seeing seems to be about $S^{[n]}$ as a moduli space of sheaves, not about Hilbert schemes \`a la Nakajima \cite{nakajima} and Grojnowski \cite{grojnowski}: we do not get a family of $\P^n$-functors $D^b(S^{[m]}) \to D^b(S^{[m+n]})$, nor do we get anything if the surface $S$ is not a K3.  We propose the following:
\begin{conj*}
Let $\M$ be a $2n$-dimensional fine moduli space of stable sheaves on a K3 surface $S$, let $F\colon D^b(S) \to D^b(\M)$ be the Fourier--Mukai functor induced by the universal sheaf $\U$ on $S \times \M$, and let $R$ be the right adjoint of $F$.  Then
\[ RF = \id_S \oplus [-2] \oplus [-4] \oplus \dotsb \oplus [-2n + 2] \]
and the monad structure $RFRF \to RF$ is like multiplication in $H^*(\P^{n-1})$, so $F$ determines an autoequivalence of $D^b(\M)$.
\end{conj*}
\noindent Of course one should be willing to drop the hypothesis that $\M$ is fine and work with twisted sheaves.  It does not seem feasible to prove this directly, as we do not know enough about the sheaves $\U|_{x \times \M}$ on $\M$, where $x \in S$; but it might be proved by deformation theory.

In \S\ref{cubic_calc} we give the following non-commutative example as evidence for the conjecture.  Let $X$ be a cubic 4-fold and let $\A = \langle \O_X, \O_X(1), \O_X(2) \rangle^\perp \subset D^b(X)$ be Kuznetsov's subcategory, which should be thought of as a non-commutative K3 surface: it has the same Serre functor and Hochschild homology and cohomology as a K3 surface \cite{kuz_cubics}, but if $X$ is generic then $\A$ lacks points and line bundles.  The variety $Y$ of lines on $X$, which is a hyperk\"ahler 4-fold, is a moduli space of objects in $\A$ by \cite[\S5]{km}, and we show that an appropriate functor $F\colon \A \to Y$ satisfies $RF = \id \oplus [-2]$, hence is spherical.  The associated autoequivalence of $D^b(Y)$ is new, as it sends the structure sheaf of a point to a complex of rank 2.\footnote{If $X$ is a generic Pfaffian cubic fourfold, so by \cite{bd} there is a K3 surface $S$ with $S^{[2]} \cong Y$, then our autoequivalences of $D^b(S^{[2]})$ and $D^b(Y)$ are in fact conjugate by tensoring by a line bundle, but the calculation is too long to include here.}

\bigskip
\vspace{-4.05ex} 
\paragraph{Acknowledgements.}  I thank Richard Thomas, Will Donovan, Ed Segal, Paul Johnson, Timothy Logvinenko, and Eyal Markman for helpful discussions, Martijn Kool for introducing me to the nested Hilbert scheme, Yujiro Kawamata for permission to use an example of his in \S\ref{P-functors}, Ciaran Meachan for his careful reading of an earlier version, Andreas Krug for pointing out a mistake, and the referee for a thorough reading and many useful suggestions.  This work was partly supported by EPSRC grant no.\ EP/G06170X/1.

\bigskip
\vspace{-4.05ex}
\paragraph{Conventions.}  All pushforwards, tensor products, etc.\ are implicitly derived.  Given objects $A$ and $B$ in a triangulated category $\B$, we write $A \perp B$ to mean that $\Hom(A,B[i]) = 0$ for all $i \in \Z$, and given $\A \subset \B$ we denote its left and right orthogonals by
\begin{align*}
^\perp\!\A &= \{ B \in \B : B \perp A \text{ for all } A \in \A \} \\
\A^\perp &= \{ B \in \B : A \perp B \text{ for all } A \in \A \}.
\end{align*}

%% file: spherical.tex

\section{Spherical functors} \label{spherical}

\subsection{Definition} \label{spherical_def}

Let $X$ be an $n$-dimensional smooth complex projective variety, and recall that an object $\E \in D^b(X)$ is called \emph{$n$-spherical} if $\Ext^*(\E,\E) \cong H^*(S^n, \C)$, where $S^n$ is the $n$-dimensional sphere, and $\E \otimes \omega_X \cong \E$.  The \emph{twist} around $\E$ is the functor $T\colon D^b(X) \to D^b(X)$ sending an object $\F$ to the cone on the evaluation map
\[ \E \otimes \RHom(\E, \F) \to \F. \]
This definition is slightly sloppy, since cones are not functorial, but by now the remedy is well-known: one can work with a dg-enhancement, or with Fourier--Mukai kernels.  We prefer the latter, so what we really mean is that $T$ is induced by the object
\[ \cone(\E^* \boxtimes \E \to \O_\Delta) \in D^b(X \times X). \]
Seidel and Thomas \cite{st} showed that $T$ is an equivalence.
\pagebreak 

Now an object of $D^b(X)$ is the same as a functor $D^b(\text{point}) \to D^b(X)$, so following Rouquier \cite{rouquier} and Anno \cite{anno}, we consider any exact functor $F\colon \A \to \B$ between triangulated categories, with left and right adjoints $L, R\colon \B \to \A$.  We define the \emph{twist} $T$ to be the cone on the counit $FR \xrightarrow\epsilon 1$ of the adjunction, so there is an exact triangle
\begin{equation} \label{counit}
FR \xrightarrow\epsilon \id_\B \to T,
\end{equation}
and the \emph{cotwist} $C$ to be the cone on the unit:
\begin{equation} \label{unit}
\id_\A \xrightarrow\eta RF \to C.
\end{equation}
(Of course we need to be in a situation where these cones make sense; we will return to this point in a moment.)  We say that $F$ is \emph{spherical} if $C$ is an equivalence and $R \cong CL$.\footnote{Rouquier requires the triangle \eqref{unit} to be split, but we do not.  Both he and Anno require a certain natural map $R \to CL$ to be an isomorphism, but this is difficult to check in practice, and in our proof of Theorem \ref{anno_thm} below we will see that any isomorphism $R \cong CL$ will do.}  If $\A$ and $\B$ have Serre functors $S_\A$ and $S_\B$, then the latter condition is equivalent to $S_\B F C \cong F S_\A$.  If $F$ is spherical then $T$ is an equivalence.

To see how this reduces to Seidel and Thomas's construction when $\A = D^b(\text{point})$, suppose that $\E \in D^b(X)$ is a spherical object and let
\[ F = \E \otimes -\colon D^b(\text{point}) \to D^b(X). \]
Then $R = \RHom(\E, -)$, so $T$ is exactly the twist defined before, and
\[ RF = \RHom(\E, \E) \otimes - = \id \oplus [-n], \]
so the cotwist $C$ is the shift $[-n]$, which is indeed an equivalence.  The condition $S_X F C \cong F S_\text{point}$ is just $\omega_X[n] \otimes \E[-n] \cong \E$.

Let us say a word about the cones \eqref{counit} and \eqref{unit}.  If $\A$ and $\B$ are derived categories of sheaves or twisted sheaves on smooth projective varieties or compact complex manifolds and $F$ is induced by a Fourier--Mukai kernel, then $R$, $RF$, and $FR$ are induced by kernels as well, the unit and counit are induced by maps of kernels, and the standard compatibilities among units and counits hold at the level of kernels \cite[Appendix]{cw}.  The same is true if $\A$ and $\B$ are admissible subcategories of these, because the projection functors are induced by kernels \cite{kuz_basechange}.\footnote{Recall that a full subcategory $\A \subset D^b(X)$ is called \emph{admissible} if the inclusion $I$ has left and right adjoints $I^l$ and $I^r$.  The main examples are the image of a fully faithful Fourier--Mukai functor and the orthogonal to an exceptional collection.  The projection functors are $I I^l$ and $I I^r$.}  It is also possible to do business with derived categories of non-compact and singular varieties if one says ``proper'' and ``perfect'' at the right moments, or with more general schemes \cite{al1}.  Rouquier's interest was in constructible sheaves.

\subsection{Examples} \label{spherical_examples}
Spherical functors unify the following special cases:
\begin{enumerate}

\item Seidel and Thomas's spherical objects, as we have discussed:  The main examples of these are a line bundle on a Calabi--Yau manifold and the structure sheaf of a $(-2)$-curve in a surface, e.g.\ $\P^1$ in its cotangent bundle.  Another is the structure sheaf of a $(-1,-1)$-curve in a 3-fold $X$, in which case the twist can also be described as doing Bondal and Orlov's flopping equivalence \cite{bo} twice: $D^b(X) \to D^b(X^+) \to D^b(X)$.

\item Horja's EZ-spherical objects \cite{horja}:  These are spherical functors of the form $F = i_*(\E \otimes q^*(-))$, where $i$ is an embedding, $q$ is a smooth bundle as in the diagram
\[ \xymatrix{
\,E\, \ar[d]_q \ar@{^(->}[r]^i & X \\
{\phantom,}Z,
} \]
and $\E \in D^b(E)$.  For example, we could take $q\colon E \to Z$ to be a $\P^1$-bundle, $X$ the total space of the relative cotangent bundle, $i\colon E \to X$ the zero section, and $\E = \O_E$.  Horja gives a sufficient condition on $\E$ for $F$ to be spherical; Huybrechts \cite[Rmk.~8.50]{huybrechts} claims that Horja's condition is stronger than necessary, and defines $\E$ to be EZ-spherical if and only if $F$ is spherical [\emph{ibid.}, Def.~8.43].

\item Toda's fat spherical objects \cite{toda}:  These are spherical functors of the form $F\colon D^b(\Spec A) \to D^b(X)$, where $A$ is an Artinian local $\C$-algebra.   Toda's first example generalizes the Atiyah flop example above to $(0,-2)$-curves.

Toda is able to simplify the hypothesis that the cotwist $C$ is an equivalence as follows.  Let $\E' \in D^b(\Spec A \times X)$ be the Fourier--Mukai kernel, let $\pi\colon \Spec A \times X \to X$ be the projection, let $0 \in \Spec A$ be the closed point, and let $\E = \E'|_{0 \times X}$.  Then his condition $\Ext^*_X(\pi_* \E', \E) \cong H^*(S^n, \C)$ is equivalent to $RF \O_0 \cong \O_0 \oplus \O_0[-n]$.  Because $\O_0$ generates $D^b(\Spec A)$, this shows that $C = [-n]$.

\newcounter{save}
\setcounter{save}{\value{enumi}}
\end{enumerate}

To date, most authors working with spherical functors have been interested in braid group representations.  Seidel and Thomas's original paper showed that on the minimal resolution of the $A_n$ surface singularity, the twists around the $(-2)$-curves satisfy the braid relations.  Khovanov and Thomas \cite{kt} constructed EZ-spherical functors from the cotangent bundles of some partial flag varieties to that of a complete flag variety and showed that the associated twists give a representation of the braid group, which they enriched to a representation of the ``braid cobordism'' 2-category.  Cautis and Kamnitzer \cite{ck1} considered a similar example and enriched the structure in a different direction, getting representations of $\mathfrak{sl}_2$ and of other Lie algebras in later papers.  Many other authors are also involved, including Rouquier and Anno; for a more complete history see \cite{ck2}.

Donovan \cite{will} gave an example in which the cotwist is more interesting than just a shift or a line bundle.  He considered certain tautological vector bundles $E_1$ on $\P^{d-1}$ and $E_2$ on $\Gr(2,d)$ and constructed a spherical functor $D^b(E_1) \to D^b(E_2)$ whose cotwist is, up to a shift and a line bundle, the twist around a spherical object on $E_1$.  He and Segal \cite{will&ed} extended this to a sequence of vector bundles $E_r$ on $\Gr(r,d)$ and spherical functors $D^b(E_r) \to D^b(E_{r+1})$ in which the twist of each is the cotwist of the next, again up to a shift and a line bundle.

The Hilbert scheme example in this paper differs from the braid group examples and Donovan's examples in that those are all EZ-spherical functors or nearly so, so objects in the image of $F$ are supported on a subvariety of the target, whereas our Fourier--Mukai kernel on $S \times S^{[2]}$ is supported everywhere.  To put it another way, if $x, y \in S$ are distinct points then $F \O_x$ and $F \O_y$ are orthogonal, but not by virtue of having disjoint support.

The cubic 4-fold example in \S\ref{cubic_calc} is unique in that the domain $\A$ of the spherical functor is not the derived category of a variety.
\bigskip

To these substantial examples we add the following ones, which are silly in that the twist is obviously an equivalence.
\begin{enumerate}

\setcounter{enumi}{\value{save}}
\item \label{i_upper_star}
Let $i\colon D \to X$ be the inclusion of a smooth divisor and take $F = i^*$, so $R = i_*$.  Then $RF = i_* i^* = i_* \O_D \otimes -$, so by rotating the exact triangle
\[ \O_X(-D) \to \O_X \to i_* \O_D \]
we find that $C = \O_X(-D)[1] \otimes -$, which is an equivalence.  The condition $R = CL$ holds because $L = i_! = i_*(\omega_i[-1] \otimes -)$ and $\omega_i = \omega_D \otimes i^* \omega_X^* = i^* \O_X(D)$ by the adjunction formula.  For the twist, by \cite[Cor.~11.4]{huybrechts} there is an exact triangle of functors
\[ (\O_D(-D)[1] \otimes -) \to i^* i_* \to \id_D, \]
so $T = \O_D(-D)[2] \otimes -$.

\item \label{i_lower_star}
Again let $i\colon D \to X$, but now take $F = i_*$, so $R = i^!$.  By a similar computation we have $C = \O_D(D)[-1] \otimes -$, $R = CL$, and $T = \O_X(D) \otimes -$.  This is the example in Anno's paper \cite{anno}.  It is an EZ-spherical twist (take $E = Z = D$), and can be seen as a family version of the fact that the skyscraper sheaf of a point in a curve is a spherical object.

This example and the previous one reflect Logvinenko's observation \cite{al2} that $F$ is spherical with cotwist $C$ and twist $T$ if and only if $R$ is spherical with cotwist $T^{-1}[1]$ and twist $C^{-1}[1]$.

\item \label{double_cover}
Let $p\colon \tilde X \to X$ be a double cover branched over a divisor $D \subset X$, and let $F = p^*$, so $R = p_*$.  Then
\[ RF = p_* p^* = p_* \O_{\tilde X} \otimes - = (\O_X \oplus \O_X(-\tfrac12 D)) \otimes -, \]
so $C = \O_X(-\tfrac12 D) \otimes -$ is an equivalence.  The condition $R = CL$ holds because $L = p_! = p_*(\omega_p \otimes -)$ and $\omega_p = p^* \O_X(\tfrac12 D)$.  For the twist, note that there is an exact triangle of functors
\[ p^* \O_X(-\tfrac12 D) \otimes \tau^* \to p^* p_* \to \id_{\tilde X}, \]
where $\tau\colon \tilde X \to \tilde X$ exchanges the two sheets of the cover, so $T = p^* \O_X(-\tfrac12 D)[1] \otimes \tau^*$.

\end{enumerate}

Seidel and Thomas prove several propositions on getting spherical objects from exceptional objects.  Recall that an object $\E \in D^b(X)$ is called \emph{exceptional} if $\Ext^*(\E, \E) = H^*(\text{point}, \C)$; the main examples are line bundles on Fano varieties and some homogeneous vector bundles.  Thus $\E$ is exceptional if and only if the functor $\E \otimes -\colon D^b(\text{point}) \to D^b(X)$ is fully faithful, so here we relate spherical functors to fully faithful functors.
\begin{prop} \label{restrict_to_admissible}
Let $\B$ be a triangulated category with Serre functor $S_\B$, and let $F\colon \B \to \mathcal C$ be a spherical functor with cotwist $C = S_\B[-k]$ for some integer $k$.  If $I\colon \A \to \B$ is fully faithful with left and right adjoints $I^l$ and $I^r$, then $F' := FI$ is spherical with cotwist $C' = S_\A[-k]$.
\end{prop}
\begin{proof}
Recall that $\A$ inherits a Serre functor from $\B$ by the formula $S_\A = I^r S_\B I$.  Let $L$ and $R$ be the adjoints of $F$, so $L' = I^l L$ and $R' = I^r R$ are the adjoints of $F'$.  The unit $\id_\A \to R'F'$ is the composition
\[ \id_\A \to I^r I \to I^r R F I, \]
and the first arrow is an isomorphism, so we find that
\[ C' = I^r C I = I^r S_\B I [-k] = S_\A [-k]. \]
Moreover the condition $S_{\mathcal C} F C = F S_\B$ is equivalent to $S_{\mathcal C} F = F[k]$, which implies $S_{\mathcal C} F' = F'[k]$, which is equivalent to $S_{\mathcal C} F' C' = F' S_\A$.
\end{proof}

From Proposition \ref{restrict_to_admissible} and our silly examples above, we recover the following examples of Seidel and Thomas:
\begin{enumerate}
\setcounter{enumi}{\value{save}}
\renewcommand \labelenumi {$\arabic{enumi}'.$}

\item Let $i\colon D \to X$ be the inclusion of an anticanonical hypersurface (that is, $\omega_X = \O_X(-D)$, so $D$ is Calabi--Yau) and let $\E \in D^b(X)$ be an exceptional object; then $i^* \E$ is spherical.  To spell things out, the set-up $\A \xrightarrow{I} \B \xrightarrow{F} \mathcal C$ from Proposition \ref{restrict_to_admissible} is
\[ D^b(\text{point}) \xrightarrow{\E \otimes -} D^b(X) \xrightarrow{i^*} D^b(D), \]
and in Example \ref{i_upper_star} above we saw that $i^*$ was spherical with cotwist $\O_X(-D)[1] \otimes - = S_X[-\dim D]$.

For example, take a smooth quartic in $\P^3$ or a smooth quintic in $\P^4$ and let $\E$ be a line bundle, or the tangent bundle.

\item Let $i\colon D \to X$ be a smooth hypersurface with $i^* \omega_X = \O_D$ and let $\E \in D^b(D)$ be an exceptional object; then $i_* \E$ is spherical.  Now we are looking at
\[ D^b(\text{point}) \xrightarrow{\E \otimes -} D^b(D) \xrightarrow{i_*} D^b(X), \]
and in Example \ref{i_lower_star} above we saw that $i_*$ was spherical with cotwist $\omega_D \otimes i^* \omega_X[-1] \otimes - = S_D[-\dim X]$.

For example, take a $(-2)$-curve in a surface and let $\E$ be a line bundle.

\item Let $p\colon X \to \P^2$ be a double cover branched over a smooth sextic, so $X$ is a K3 surface, and let $\E \in D^b(\P^2)$ be an exceptional object; then $p^* \E$ is spherical.  Now we are looking at
\[ D^b(\text{point}) \xrightarrow{\E \otimes -}  D^b(\P^2) \xrightarrow{p^*} D^b(X), \]
and in Example \ref{double_cover} above we saw that $p^*$ was spherical with cotwist $\O_{\P^2}(-3) \otimes - = S_{\P^2}[-2]$.

\end{enumerate}

\subsection{Splitting of \texorpdfstring{$FRF$}{FRF}} \label{FRF}

The following simple observation will be the key to describing the action of $T$ on cohomology, proving that $T$ is an equivalence, and constructing the $\P$-twist associated to a $\P$-functor in \S\ref{def_of_P}.  While the unit $\eta\colon \id_\B \to RF$ is not split in general, the map $F\eta\colon F \to FRF$ is naturally split: we have a commutative triangle
\[ \xymatrix{
F \ar@{=}[rd] \ar[r]^{F \eta}
& FRF \ar[d]^{\epsilon F} \\
& F.
} \]
In the down-to-earth case $F = \E \otimes -: D^b(\text{point}) \to D^b(X)$, we are saying that the map $\C \to \RHom(\E,\E)$ that sends 1 to the identity is not split by any natural map $\RHom(\E,\E) \to \C$,\footnote{The trace map almost works, but not when $\rank \E = 0$.} but if we tensor with $\E$ then the map $\E \to \E \otimes \RHom(\E,\E)$ is split by the evaluation map $\E \otimes \RHom(\E,\E) \to \E$.

Extend the commutative triangle above to
\[ \xymatrix{
& TF[-1] \ar[d] \ar[rd] \\
F \ar@{=}[rd] \ar[r]^{F \eta}
& FRF \ar[d]^{\epsilon F} \ar[r]
& FC \\
& F,
} \]
where the row and column are exact.  Using the octahedral axiom we see that the upper-right diagonal map is an isomorphism:
\begin{equation} \label{TFFC}
TF[-1] \xrightarrow\cong FC.
\end{equation}
Note that this is true for any $F$, spherical or not.  Similarly, by looking at $RFR$, $FLF$, and $LFL$ we get isomorphisms
\begin{align}
RT[-1] &\xrightarrow\cong CR \label{RTCR} \\
F C^l &\xrightarrow\cong T^l F[1] \label{FClTlF} \\
C^l L &\xrightarrow\cong L T^l[1], \label{ClLLTl}
\end{align}
where $T^l$ and $C^l$ are the left adjoints of $T$ and $C$.

While we are here we make one more observation, which we will need in \S\ref{ploog_proof} and \S\ref{def_of_P}.  We have just seen that $FRF$ splits as $F \oplus TF[-1]$ or $F \oplus FC$.  Thus the identity map $FRF \to FRF$ can be written as the sum of two idempotents, namely the compositions
\begin{gather*}
FRF \xrightarrow{\epsilon F} F \xrightarrow{F \eta} FRF \\
FRF \to FC \xleftarrow\cong TF[-1] \to FRF.
\end{gather*}

\subsection{Action on spanning class and cohomology} \label{spherical_coho}

If $\E$ is a spherical object on $X$ then $\{ \E \} \cup \E^\perp$ is a spanning class for $D^b(X)$; that is, an object that is left and right orthogonal to $\E$ and $\E^\perp$ is zero, although not every object can be gotten from $\E$ and $\E^\perp$ by taking cones.  The twist $T$ sends $\E$ to $\E[-n+1]$ and acts on $\E^\perp$ as the identity \cite[Ex.~8.5(ii)]{huybrechts}.  If $X$ is even-dimensional, the induced action on cohomology is a reflection, sending the Mukai vector $v(\E) \in H^*(X, \Q)$ to $-v(\E)$ and acting as the identity on its orthogonal $v(\E)^\perp$ under the Mukai pairing.

For a spherical functor $F\colon \A \to \B$, this is generalized as follows.  We replace $\E$ with the set of objects
\[ \im F = \{ FA : A \in \A \}, \]
and $\E^\perp$ with
\[ \ker R = \{ B \in \B : RB = 0 \}. \]
It is easy to see that $\ker R = (\im F)^\perp$.  I claim that $\im F \cup \ker R$ is a spanning class for $\B$.  First, if $(\im F) \perp B$ then $B \in \ker R$, and if in addition $(\ker R) \perp B$ then $B \perp B$, so $B = 0$; thus $(\im F \cup \ker R)^\perp = 0$.  On the other side we see that $^\perp(\im F) = \ker L$, but since $R \cong CL$ and $C$ is an equivalence we see that $\ker L = \ker R$; thus by a similar argument we find that $^\perp(\im F \cup \ker R) = 0$.

The twist $T$ acts as the identity on $\ker R$, for if $B \in \ker R$ then the first term in the exact triangle
\[ FRB \to B \to TB \]
vanishes, so $TB = B$.  It acts on $\im F$ by
\[ TFA = FCA[1] \]
as we saw in the previous section.  In particular, if $C = [-n]$ then $T$ acts on $\im F$ by $[-n+1]$.
\bigskip

We are now in a position to describe how the autoequivalence discussed at the beginning of the introduction acts on cohomology.  Recall that $S$ is a K3 surface and $F\colon D^b(S) \to D^b(S^{[2]})$ is induced by the universal ideal sheaf, and we will prove later that $RF = \id \oplus [-2]$.  The induced map
\[ F^h\colon H^*(S,\Q) \to H^*(S^{[2]}, \Q) \]
is injective, since $R^h F^h = (RF)^h$ is multiplication by 2, and
\[ T^h\colon H^*(S^{[2]},\Q) \to H^*(S^{[2]}, \Q) \]
acts as multiplication by $-1$ on $\im F^h$ and as the identity on its orthogonal $(\im F^h)^\perp$ under the Mukai pairing.

We can also show that our twist does not come from any known spherical twist on $S$ via Ploog's map
\[ \varphi\colon \Aut(D^b(S)) \hookrightarrow \Aut(D^b(S^{[2]})). \]
This map uses the Bridgeland--King--Reid--Haiman equivalence $D^b(S^{[2]}) \cong D^b([S^2/\mathfrak S_2])$, where the latter is the derived category of the quotient stack, or equivalently the $\mathfrak S_2$-equivariant derived category of $S^2$.  Suppose that $\E \in D^b(S)$ is a spherical object and $\F \in \E^\perp$.\footnote{One does not know whether such an $\F$ exists for an arbitrary spherical object $\E$, but it does exist in all known examples.  If $\E$ is a line bundle, take $\F = \E \otimes \I_x^* \otimes \I_y$, where $x, y \in S$ are distinct points.  If $\E$ is the structure sheaf of a $(-2)$-curve, take $\F = \O_x$ for some point $x$ not on the curve.  For a construction of $\F$ when $\E$ is an arbitrary stable vector bundle, see \cite[Example~1.24]{ploog_thesis}.}  The spherical twist $T_\E$ shifts $\E$ by $-1$ and fixes $\F$.  Consider the objects 
\begin{equation} \label{ploog_objects}
\left. \begin{array}{c}
\E \boxtimes \E \\
(\E \boxtimes \F) \oplus (\F \boxtimes \E) \\
\F \boxtimes \F
\end{array} \right\}
\in D^b([S^2/\mathfrak S_2]).
\end{equation}
Then $\varphi(T_\E)$ shifts the first by $-2$, shifts the second by $-1$, and fixes the third.  On the other hand, our twist $T_F$ shifts $\im F$ by $-1$ and fixes $\ker R$.\footnote{We cannot rule out the possibility that $\ker R = 0$, but this does not affect the argument that follows.}  Now we need the following:
\begin{lem}
Suppose that $X$ is a smooth variety, $A, B \in D^b(X)$, and $T$ is an autoequivalence of $D^b(X)$ with $TA = A[i]$ and $TB = B[j]$ for some $i \ne j \in \Z$.  Then $A \perp B$ and $B \perp A$.
\end{lem}
\begin{proof}
For all $k, m \in \Z$ we have
\begin{align*}
\Hom(A,B[k])
&= \Hom(T^m A, T^m B[k]) \\
&= \Hom(A[mi], B[k+mj]) \\
&= \Hom(A, B[k+m(j-i)]),
\end{align*}
which vanishes for $|m| \gg 0$ because $X$ is smooth.  Similarly, $\Hom(B,A[k]) = 0$ for all $k$.
\end{proof}
\noindent Thus if $T_F$ were $\varphi(T_\E)$ or a shift of it, then one of the objects in \eqref{ploog_objects} would be orthogonal to the spanning class $\im F \cup \ker R$, hence would be zero, which gives a contradiction.

In the introduction we observed that $T_F$ is not generated by shifts, line bundles, automorphisms of $S^{[2]}$, or twists around $\P^2$-objects, because these all preserve rank (up to sign) while $T_F$ sends the structure sheaf of a point to a rank-2 sheaf shifted by 2.  Another known autoequivalence of $D^b(S^{[2]})$ is the following EZ-spherical twist \cite[Example 8.49(iv)]{huybrechts}:  Consider
\[ E = \{ \xi \in S^{[2]} : \text{$\supp \xi$ is a single point} \}, \]
which is the exceptional divisor of the Hilbert--Chow morphism $S^{[2]} \to S^{(2)}$.  It is a $\P^1$-bundle over $S$ -- the projectivization of the tangent bundle, in fact.  Let $q\colon E \to S$ be the $\P^1$-bundle, and let $i\colon E \to S^{[2]}$ be the inclusion.  Then from Proposition \ref{restrict_to_admissible} and the examples in \S\ref{spherical_examples} we easily check that $i_* q^*\colon D^b(S) \to D^b(S^{[2]})$ is spherical with cotwist $[-2]$.  But EZ-spherical twists preserve rank as well: if $\F \in D^b(S^{[2]})$ then from the triangle
\[ i_* q^* q_* i^! \F \to \F \to T_{i_* q^*} \F \]
we see that
\begin{align*}
\rank T_{i_* q^*} \F &= \rank \F - \rank i_* q^* q_* i^! \F \\
&= \rank\F - 0.
\end{align*}

Of course one would like to know whether $T_F$ is in the subgroup generated by these rank-preserving equivalences and the image of Ploog's map $\varphi$, but this question is too difficult to settle at present.

\pagebreak
\subsection{Proof of equivalence} \label{ploog_proof}
We conclude with an alternate proof that $T$ is an equivalence, following Ploog \cite[Thm.~1.27]{ploog_thesis}.
\begin{thm}[Rouquier, Anno] \label{anno_thm}
If $F\colon \A \to \B$ is spherical then the twist $T$ is an equivalence.
\end{thm}
\begin{proof}
In the previous section we saw that $\im F \cup \ker R$ is a spanning class for $\B$.  By \cite[Prop.~1.49]{huybrechts}, we can show that $T$ is fully faithful by showing that the natural map
\[ \Hom(B,B'[i]) \to \Hom(TB, TB'[i]) \]
is an isomorphism for all $B, B' \in (\im F \cup \ker R)$ and all $i \in \Z$.  Since $\im F \cup \ker R$ is closed under shifts, we need only consider $i=0$.

We check this in four cases.  First, if $B, B' \in \ker R$ then $TB = B$ and $TB' = B'$, as we saw in the previous section, so $\Hom(TB, TB') = \Hom(B,B')$.  Next, if $FA \in \im F$ and $B \in \ker R = \ker L$ then
\begin{align*}
\Hom(TFA, TB)
&= \Hom(FCA[1], B) \\
&= \Hom(CA[1], RB) \\
&= 0 \\
&= \Hom(FA, B) \\[1em]
\Hom(TB, TFA)
&= \Hom(B, FCA[1]) \\
&= \Hom(LB, CA[1]) \\
&= 0 \\
&= \Hom(B, FA).
\end{align*}
Last, if $FA, FA' \in \im F$ then
\begin{equation} \label{chain_of_eqs}
\begin{split}
\Hom(TFA, TFA')
&= \Hom(T^lTFA, FA') \\
&= \Hom(T^lFCA[1], FA') \\
&= \Hom(FC^lCA, FA') \\
&= \Hom(FA, FA')
\end{split}
\end{equation}
where in the second line we have use \eqref{TFFC}, in the third we have used \eqref{FClTlF}, and in the last we have $C^l C = \id_\A$ because $C$ is an equivalence.  But this is not quite enough: we must show that
\[ T^l T F \xrightarrow{\epsilon F} F \]
is an isomorphism.  The chain of equalities \eqref{chain_of_eqs} suggests showing that it equals the composition
\[ T^l T F \cong T^l F C[1] \cong F C^l C \xrightarrow{F \epsilon} F. \]
This is terribly boring, and we prove it as a separate lemma below; in fact they are the same up to a sign, which is good enough.

Now $T$ is fully faithful, so by \cite[Ex.~1.51]{huybrechts} we can show that it is an equivalence by showing that $\ker T^l = 0$.  If $B \in \ker T^l$ then $C^l L B = L T^l B[1] = 0$, but $C^l$ is an equivalence, so $LB = 0$.  Take left adjoints of \eqref{counit} to get an exact triangle
\[ T^l \to \id_\B \to FL, \]
from which we see that if $T^l B = 0$ then $B = 0$.
\end{proof}

\begin{lem}
Let $F\colon \A \to \B$ be a functor, not necessarily spherical, with left and right adjoints $L$ and $R$, let $T$ and $C$ be the twist and cotwist as in \S\ref{spherical_def}, and let $T^l$ and $C^l$ be their left adjoints.  Then the compositions
\begin{align*}
T^l F C[1] &\cong T^l T F \xrightarrow{\epsilon F} F \\
T^l F C[1] &\cong F C^l C \xrightarrow{F \epsilon} F,
\end{align*}
where the isomorphisms $\cong$ are as in \eqref{TFFC} and \eqref{FClTlF}, are equal up to a sign.
\end{lem}
\begin{proof}
First note that for any two functors $\Phi, \Psi\colon \mathcal C \to \mathcal D$ with left adjoints $\Phi^l, \Psi^l\colon \mathcal D \to \mathcal C$, a natural transformation $\tau\colon \Phi \to \Psi$ determines a natural transformation $\tau^l\colon \Phi^l \to \Psi^l$, and the diagram
\[ \xymatrix{
\Psi^l \Phi \ar[d]_{\tau^l \Phi} \ar[r]^{\Psi^l \tau} & \Psi^l \Psi \ar[d]^\epsilon \\
\Phi^l \Phi \ar[r]_\epsilon & \id_\mathcal C
} \]
commutes.  This is clear if $\mathcal C$ and $\mathcal D$ are the derived categories of smooth compact spaces, $\Phi$ and $\Psi$ are Fourier--Mukai functors, and $\tau$ is induced by a map of kernels, because $\Phi^l$ and $\Psi^l$ are induced by the dual kernels (tensored with the shift of a line bundle).  But it is true in any category: $\tau^l$ is the composition $\Psi^l \xrightarrow{\Psi^l \eta}
\Psi^l \Phi \Phi^l \xrightarrow{\Psi^l \tau \Phi^l} \Psi^l \Psi \Phi^l \xrightarrow{\epsilon \Phi^l} \Phi^l$, and it is easy to check that the diagram above commutes.

Consider the diagram\footnote{This diagram contains no exact triangles.}
\[ \xymatrix{
& F C^l RF \ar[d] \ar[r] & F C^l C \ar[d] \ar[rrdd]|-{F \epsilon} \\
FLTF[-1] \ar[d] \ar[r] & FLFRF \ar[d] \ar[rddd]|-{\epsilon F} \ar[rrrd]|-{F \epsilon} \ar[r] & FLFC \ar[d] \\
T^l T F \ar[rrdd]|-{\epsilon F} \ar[r] & T^l FRF[1] \ar[r] & T^l FC[1] & & F\\
\\
& & F.
} \]
The two ``kites'' (involving $\epsilon F$ and $F \epsilon$) commute by the preceding discussion, taking $\tau\colon \Phi \to \Psi$ to be $T[-1] \to FR$ or $RF \to C$.  The three squares obviously commute.  The two horizontal compositions and the two vertical compositions are isomorphisms as we saw in \eqref{TFFC} and \eqref{FClTlF}.  Thus the composition $FLFRF \to T^lFC[1]$ is an epimorphism, so to prove the lemma it is enough to show that the compositions
\begin{align}
FLFRF &\to T^l F C[1] \xleftarrow\cong T^l T F \xrightarrow{\epsilon F} F \label{comp1} \\
FLFRF &\to T^l F C[1] \xleftarrow\cong F C^l C \xrightarrow{F \epsilon} F \label{comp2}
\end{align}
are equal up to a sign.

With reference to the big diagram above, we can rewrite \eqref{comp1} in the following steps:
\begin{gather}
FLFRF \to FLFC \to T^l F C[1] \xleftarrow\cong T^l T F \xrightarrow{\epsilon F} F \notag \\
FLFRF \to FLFC \xleftarrow\cong FLTF[-1] \to T^l T F \xrightarrow{\epsilon F} F \notag \\
FLFRF \to FLFC \xleftarrow\cong FLTF[-1] \to FLFRF \xrightarrow{\epsilon F} F. \label{step_3}
\end{gather}
In \S\ref{FRF} we saw that the idempotent
\[ FLFRF \to FLFC \xleftarrow\cong FLTF[-1] \to FLFRF, \]
which is the first three steps of \eqref{step_3}, equals the identity minus
\begin{equation} \label{step_3a}
FLFRF \xrightarrow{FL\epsilon F} FLF \xrightarrow{FLF\eta} FLFRF.
\end{equation}
The map $FLFRF \xrightarrow{\epsilon F} F$, which is the last step of \eqref{step_3}, can be factored as
\begin{equation} \label{step_3b}
FLFRF \xrightarrow{F\epsilon RF} FRF \xrightarrow{\epsilon F} F.
\end{equation}
The composition of \eqref{step_3a} and \eqref{step_3b} is just $FLFRF \xrightarrow{F \epsilon} F$, as we see from the following diagram:
\[ \xymatrix{
& & FLFRF \ar[rd]^{F\epsilon RF} \\
FLFRF \ar[r]^{FL\epsilon F} & FLF \ar[ru]^{FLF\eta} \ar[rd]_{F \epsilon} & & FRF \ar[r]^{\epsilon F} & F. \\
& & F \ar[ru]^{F \eta} \ar@{=}@/_/[rru]
} \]
We conclude that \eqref{step_3}, and hence \eqref{comp1}, equals $FLFRF \xrightarrow{\epsilon F - F \epsilon} F$.

Similarly we find that \eqref{comp2} equals $FLFRF \xrightarrow{F \epsilon - \epsilon F} F$.
\end{proof}

%% file: hilb_calc.tex

\pagebreak
\section{Hilbert scheme calculation} \label{hilb_calc}

In this section we prove the following:
\begin{thm} \label{hilb_calc_thm}
Let $S$ be a complex projective K3 surface, $S^{[n]}$ its Hilbert scheme of length-$n$ subschemes, $Z = Z_n \subset S \times S^{[n]}$ the universal subscheme, $F\colon D^b(S) \to D^b(S^{[n]})$ the functor induced by the ideal sheaf $\I_Z$, and $R$ the right adjoint of $F$. 
\begin{enumerate}
\renewcommand \labelenumi {(\alph{enumi})}
\item There is an isomorphism
\[ R F \cong \id_S \oplus [-2] \oplus [-4] \oplus \dotsb \oplus [-2n+2].  \]
\item This isomorphism can be chosen so that the map
\[ RF[-2] \hookrightarrow RFRF \xrightarrow{R \epsilon F} RF, \]
when written in components
\[ [-2] \oplus [-4] \oplus \dotsb \oplus [-2n] \longrightarrow \id_S \oplus [-2] \oplus \dotsb \oplus [-2n+2], \]
is of the form
\[ \begin{pmatrix}
0  \\
1 & 0 \\
  & 1 & 0 \\
& & \ddots & \ddots \\
& & & 1 & 0 & * \\
& & & & 1 & *
\end{pmatrix}. \]
\end{enumerate}
\end{thm}

In addition to $F$, we will consider the functors $F', F''\colon D^b(S) \to D^b(S^{[n]})$ induced by $\O_{S \times S^{[n]}}$ and $\O_Z$ respectively, and their right adjoints $R'$ and $R''$.  We have exact triangles of functors
\begin{align*}
F &\to F' \to F'' \\
R'' &\to R' \to R.
\end{align*}
In \S\ref{nested} we give an exposition of the ``nested Hilbert scheme'' which will be central to our computations.  In \S\S\ref{easy}--\ref{cancel} we compute $R'F'$, $R'F''$, $R''F'$, and $R''F''$, and enough information about the maps between them to determine $RF$ through some long exact sequences.  In \S\ref{monad} we prove statement (b) about the monad structure of $RF$.

\subsection{Nested Hilbert schemes} \label{nested}

The nested Hilbert scheme is
\[ S^{[n-1,n]} = \{ (\zeta, \xi) \in S^{[n-1]} \times S^{[n]} : \zeta \subset \xi \}. \]
Like $S^{[n]}$, it is $2n$-dimensional and smooth \cite{tik}.  We give a quick tour of its geometry, following Ellingsrud and Str{\o}mme \cite{es}.  This discussion is valid for any smooth surface.

For motivation, recall that $S^{[2]}$ has a very simple construction: let $\Delta \subset S \times S$ be the diagonal; then the involution of $S \times S$ lifts to $\Bl_\Delta(S \times S)$, fixing the exceptional divisor $E$, and the quotient is $S^{[2]}$.  We summarize this in the diagram
\[ \xymatrix{
E \ar@{}[r]|-*{\subset} \ar[d]
& \Bl_\Delta(S \times S) \ar[r]^{\qquad g} \ar[d]^\gamma
& S^{[2]} \\
\Delta \ar@{}[r]|-*{\subset}
& S \times S.
} \]
The map $\pi_1 \gamma \times g\colon \Bl_\Delta(S \times S) \to S \times S^{[2]}$ is an embedding, and its image is the universal subscheme $Z_2$.

For $n>2$, the picture will be
\[ \xymatrix{
E \ar@{}[r]|-*{\subset} \ar[d] 
& S^{[n-1,n]} \ar[r]^{\qquad g} \ar[d]^{\gamma = q \times f}
& S^{[n]} \\
Z_{n-1} \ar@{}[r]|-*{\subset}
& S \times S^{[n-1]},
} \]
where $f\colon S^{[n-1,n]} \to S^{[n-1]}$ and $g\colon S^{[n-1,n]} \to S^{[n]}$ are the obvious maps and $q\colon S^{[n-1,n]} \to S$ sends a pair $\zeta \subset \xi$ to the point where they differ, that is, where the kernel of $\O_\xi \to \O_\zeta$ is supported, which we will call $\xi \setminus \zeta$.\footnote{But note that there is no similar map $S^{[n-m,n]} \to S^{[m]}$ for $m > 1$, because the kernel of $\O_\xi \twoheadrightarrow \O_\zeta$ need not be a quotient of $\O_S$.}

Let $\phi = q \times g\colon S^{[n-1,n]} \to S \times S^{[n]}$.  For any $(\zeta,\xi) \in S^{[n-1,n]}$ we have an exact sequence
\[ 0 \to \O_{\xi\setminus\zeta} \to \O_\xi \to \O_\zeta \to 0 \]
so we see that the fiber of $\phi$ over $(x, \xi)$ is $\P\Hom(\O_x, \O_\xi)^*$.\footnote{In this section only, we use Grothendieck's convention that $\P$ is the projective space of 1-dimensional quotients.  The reason will be clear in the next paragraph.}  Thus the image of $\phi$ is $Z_n$, and $\phi$ is an isomorphism over the set of $(x, \xi) \in Z_n$ where the length of $\xi$ at $x$ is 1, so $\phi$ is a resolution of singularities for $Z_n$.  Since the fibers of $\phi$ are projective spaces, we have
\[ \phi_* \O_{S^{[n-1,n]}} = \O_{Z_n}, \]
so $Z_n$ has rational singularities.

Next let $\gamma = q \times f\colon S^{[n-1,n]} \to S \times S^{[n-1]}$.  For any $(\zeta,\xi) \in S^{[n-1,n]}$ we have an exact sequence
\[ 0 \to \I_\xi \to \I_\zeta \to \O_{\xi\setminus\zeta} \to 0 \]
so we see that the fiber of $\gamma$ over $(x, \zeta)$ is $\P\Hom(\I_\zeta, \O_x)^* = \P(\I_\zeta|_x)$, so $S^{[n-1,n]}$ is isomorphic to the projectivization\footnote{For the reader who is uncomfortable with projectivizing sheaves that are not vector bundles, I recommend \cite[pp.~103, 115, and 170--171]{eh}.}
\[ \P \I_{Z_{n-1}} = \Proj(\O_{S \times S^{[n-1]}} \oplus \I_{Z_{n-1}} \oplus \Sym^2 \I_{Z_{n-1}} \oplus \dotsb). \]
The blowup
\[ \Bl_{Z_{n-1}}(S \times S^{[n-1]}) = \Proj(\O_{S \times S^{[n-1]}} \oplus \I_{Z_{n-1}} \oplus \I_{Z_{n-1}}^2 \oplus \dotsb) \]
naturally embeds into $\P \I_{Z_{n-1}} \cong S^{[n-1,n]}$, and since the latter is smooth, hence irreducible, the embedding is an isomorphism.  Note that the rational map $g \circ \gamma^{-1}\colon \xymatrix{S \times S^{[n-1]} \ar@{-->}[r] & S^{[n]}}$ just sends a pair $(x,\zeta) \notin Z_{n-1}$ to $x \cup \zeta$.

Now $Z_{n-1}$ is singular for $n>3$, and it is perhaps strange to blow up a smooth variety along a singular center and end up with a smooth variety.  But $\gamma$ behaves in many ways like a blowup along a smooth center:
\begin{prop} \label{gamma_prop} {\ }
\begin{enumerate}
\renewcommand \labelenumi {(\alph{enumi})}
\item $\gamma_* \O_{S^{[n-1,n]}} = \O_{S \times S^{[n-1]}}$.
\item The exceptional divisor
\[ E = \gamma^{-1}(Z_{n-1}) = \{ (\zeta, \xi) \in S^{[n-1,n]} : (\xi \setminus \zeta) \in \zeta \} \]
is irreducible.
\item The relative canonical bundle $\omega_\gamma = \O_{S^{[n-1,n]}}(E)$.
\item $\gamma_* \O_E(E) = 0$.
\end{enumerate}
\end{prop}
\begin{proof} \ 
\begin{enumerate}
\renewcommand \labelenumi {(\alph{enumi})}
\item The fibers of $\gamma$ are projective spaces.
\item This is proved in \cite[\S3]{es}.
\item Let $Z_{n-1}'$ be the singular locus of $Z_{n-1}$, and let $E' = \gamma^{-1}(Z_{n-1}')$.  Then $\gamma\colon (S^{[n-1,n]} \setminus E') \to (S \times S^{[n-1]} \setminus Z_{n-1}')$ is a blow-up along a smooth center of codimension 2, so the line bundles $\omega_\gamma$ and $\O_{S^{[n-1,n]}}(E)$ agree away from $E'$.  But $E'$ is a proper subset of the irreducible divisor $E$, so it has codimension at least 2 in $S^{[n-1,n]}$, so the claim follows by Hartogs' theorem.
\item Take the exact sequence
\[ 0 \to \O_{S^{[n-1,n]}} \to \O_{S^{[n-1,n]}}(E) \to \O_E(E) \to 0 \]
and apply $\gamma_*$ to get an exact triangle
\[ \O_{S \times S^{[n-1]}} \to \O_{S \times S^{[n-1]}} \to \gamma_* \O_E(E). \]
The first map is an isomorphism away from $\gamma(E) = Z_{n-1}$, which has codimension 2, so it is an isomorphism. \qedhere
\end{enumerate}
\end{proof}

We conclude with the following fact:
\begin{prop} The map $q$ is a submersion. \end{prop}
\begin{proof}
This can be proved by working directly with the tangent spaces, but the proof is messy.  Instead we give a quick transcendental proof.  By Sard's theorem, $q$ is a submersion over almost all $x \in S$.  If $S = \C^2$, this implies that $q$ is a submersion everywhere by translation.  Now for any smooth surface $S$, let $(\zeta,\xi) \in S^{[n-1,n]}$ and let $U \subset S$ be an analytic neighborhood of $\supp \xi$ isomorphic to an open set in $\C^2$, possibly disconnected.  Then $U^{[n-1,n]}$ is a neighborhood of $(\zeta,\xi)$, and we have
\[ \xymatrix{
U^{[n-1,n]} \ar@{^(->}[r] \ar[d]_q & (\C^2)^{[n-1,n]} \ar[d]^q \\
U \ar@{^(->}[r] & \C^2.
} \]
The horizontal maps are embeddings of open sets, and we have just seen that the right-hand $q$ is a submersion.
\end{proof}

\subsection{\texorpdfstring{$\mathbf{R'F', R'F'', R''F'}$}{R'F', R'F'{'},R'{'}F'}, and the maps between them} \label{easy}

We return to the setting of Theorem \ref{hilb_calc_thm}, so $S$ is a K3 surface and $Z$ is short for $Z_n \subset S \times S^{[n]}$.

\paragraph{$\mathbf{R'F'}$.}  It will be convenient to use the same name to refer to a functor and the kernel that induces it; thus
\begin{align*}
F' &= \O_{S \times S^{[n]}} \\
R' &= \O_{S^{[n]} \times S}[2] \\
R'F' &= \pi_{SS*} \O_{S \times S^{[n]} \times S} [2]
\end{align*}
where $\pi_{SS}$ is the projection $S \times S^{[n]} \times S \to S \times S$.  Then we see that
\begin{align*}
R'F' &= \O_{S \times S} \otimes \RGamma(\O_{S^{[n]}})[2] \\
&= \O_{S\times S}[2] \oplus \O_{S \times S} \oplus \O_{S\times S}[-2] \oplus \dotsb \oplus \O_{S\times S}[-2n+2].
\end{align*}

\paragraph{$\mathbf{R'F''}$.}  Next we have $F'' = \O_Z$, so $R'F'' = \pi_{SS*} \O_{S \times Z}[2]$.  Consider the diagram
\begin{equation} \label{pushing_OZ}
\begin{split}
\xymatrix{
S^{[n-1,n]} \ar[d]_{\gamma = f \times q} \ar[r]^{\phi = g \times q}
& S^{[n]} \times S \ar[dd]^{\pi_S} \\
S^{[n-1]} \times S \ar[rd]_{\varpi_S} \\
& S.
} \end{split}
\end{equation}
Then we have
\begin{align*}
R'F'' &= \O_S \boxtimes \pi_{S*} \O_Z[2] \\
&= \O_S \boxtimes \pi_{S*} \phi_* \O_{S^{[n-1,n]}}[2] \\
&= \O_S \boxtimes \varpi_{S*} \gamma_* \O_{S^{[n-1,n]}}[2] \\
&= \O_S \boxtimes \varpi_{S*} \O_{S^{[n-1]} \times S}[2] \\[.5\baselineskip]
&= \O_{S \times S} \otimes \RGamma(\O_{S^{[n-1]}})[2] \\
&= \O_{S \times S}[2] \oplus \O_{S \times S} \oplus \O_{S \times S}[-2] \oplus \dotsb \oplus \O_{S \times S}[-2n+4].
\end{align*}

\paragraph{$\mathbf{R'F}$.}  Next I claim that the map $R'F' \to R'F''$ induces an isomorphism on $\H^i$ for $i < 2n-2$, so
\[ R'F = \O_{S\times S}[-2n+2]. \]
This amounts to claiming that in the diagram \eqref{pushing_OZ}, the restriction map $\O_{S^{[n]} \times S} \to \phi_* \O_{S^{[n-1,n]}}$ induces an isomorphism on $R^i \pi_{S*}$ for $i<2n$.  We check this fiberwise.  Since $q$ is a submersion, its fibers are smooth.  Over a point $x \in S$, the fiber of \eqref{pushing_OZ} is
\[ \xymatrix{
q^{-1}(x) \ar[d]_f \ar[r]^g
& S^{[n]} \ar[dd] \\
S^{[n-1]} \ar[rd] \\
& x.
} \]
Now we want to show that $g^*\colon H^i(\O_{S^{[n]}}) \to H^i(\O_{q^{-1}(x)})$ is an isomorphism for $i<2n$.  Let $\sigma$ be a non-vanishing holomorphic 2-form on $S$, and let $\sigma_{n-1}$ and $\sigma_n$ be the induced holomorphic 2-forms on $S^{[n-1]}$ and $S^{[n]}$ constructed by Beauville \cite[Prop.\ 5]{beauville}.  From his construction it is easy to check that on $S^{[n-1,n]}$ we have $g^* \sigma_n = q^* \sigma + f^* \sigma_{n-1}$.  Thus the generator $\bar\sigma_n^j$ of $H^{2j}(\O_{S^{[n]}})$ maps to $f^* \bar\sigma_{n-1}^j \in H^{2j}(\O_{q^{-1}(x)})$.  But since $f_* \O_{q^{-1}(x)} = \O_{S^{[n-1]}}$, the map $f^*\colon H^i(\O_{S^{[n-1]}}) \to H^i(\O_{q^{-1}(x)})$ is an isomorphism, so $f^* \bar\sigma_{n-1}^j$ generates $H^{2j}(\O_{q^{-1}(x)})$ for $j<n$, as desired.

\paragraph{$\mathbf{R''F'}$ and $\mathbf{RF'}$.}  By duality we have
\begin{align*}
R''F' &= \O_{S \times S} \otimes \RGamma(\O_{S^{[n-1]}}) \\
&= \O_{S \times S} \oplus \O_{S \times S}[-2] \oplus \dotsb \oplus \O_{S \times S}[-2n+2],
\end{align*}
and the map $R''F' \to R'F'$ induces an isomorphism on $\H^i$ for $i > -2$, so
\[ RF' = \O_{S \times S}[2]. \]

\pagebreak
\subsection{Main calculation: \texorpdfstring{$\mathbf{R''F''}$}{R'{'}F'{'}}} \label{hard}

In this section we show that
\begin{align}
R''F'' &= (\O_\Delta \otimes \RGamma(\O_{S^{[n-1]}})) \oplus (\O_{S \times S} \otimes \RGamma(\O_{S^{[n-2]}})) \label{R''F''} \\[2pt]
&= (\O_\Delta \oplus \O_\Delta[-2] \oplus \dotsb \oplus \O_\Delta[-2n+4] \oplus \O_\Delta[-2n+2]) \notag \\
&\qquad\qquad\qquad \oplus (\O_{S \times S} \oplus \O_{S\times S}[-2] \oplus \dotsb \oplus \O_{S\times S}[-2n+4]). \notag
\end{align}
The essential reason is as follows.  We have
\[ R''F'' = \pi_{SS*}(\O_{Z\times S} \otimes \O_{S\times Z}^*[2]), \]
where the tensor product is taken on $S \times S^{[n]} \times S$, and $\O_{Z\times S} \otimes \O_{S\times Z}^*$ is supported on $(Z \times S) \cap (S \times Z) \cong Z \times_{S^{[n]}} Z$, which has two irreducible components: the diagonal $Z$, and the rest, which is birational to $S \times S \times S^{[n-2]}$.  These two components are responsible for the two summands of \eqref{R''F''}.  We mention this now for fear that it will be obscured in the computation that follows.

To carry out the computation, replace $Z \times_{S^{[n]}} Z$ with the following partial desingularization:
\begin{align*}
X &:= Z \times_{S^{[n]}} S^{[n-1,n]} \\
&= \{ (x, \zeta, \xi) \in S \times S^{[n-1,n]} : x \in \xi \}.
\end{align*}
The diagram
\[ \xymatrix{
X \ar[d]_{\tilde\phi'} \ar[rr]^-{\tilde\imath}
&& S \times S^{[n-1,n]} \ar[d]^{\phi' := 1 \times g \times q} \\
Z \times S \ar[rr]_-i
&& S \times S^{[n]} \times S
} \]
is Cartesian, and we have
\begin{align*}
\O_{Z\times S} \otimes \O_{S\times Z}^*[2]
&= \H om(\phi'_* \O_{S \times S^{[n-1,n]}}, i_* \O_{Z \times S})[2] \\
&= \phi'_*(\phi'^* i_* \O_{Z \times S} \otimes \pi_2^* \O(E)) \\
&= \phi'_*(\tilde\imath_* \O_X \otimes \pi_2^* \O(E)).
\end{align*}
In the second line, $\pi_2$ is the projection $S \times S^{[n-1,n]} \to S^{[n-1,n]}$, and we have used Grothendieck duality: from Proposition \ref{gamma_prop}(c) we know that $\O_{S^{[n-1,n]}}(E)$ is the relative canonical bundle of $\gamma\colon S^{[n-1,n]} \to S \times S^{[n-1]}$, hence is the canonical bundle of $S^{[n-1,n]}$, so $\pi_2^* \O(E)$ is the canonical bundle of $S \times S^{[n-1,n]}$, hence is the relative canonical bundle of $\phi'$.  In the third line we have used the base change criterion in Appendix \ref{base_change}, which requires that every irredicible component of $X$ have dimension $2n$; to see that this is true, observe that $Z$ is flat and finite over $S^{[n]}$, so $X$ flat and finite over $S^{[n-1,n]}$, and since the latter is smooth, $X$ is Cohen--Macaulay, hence equidimensional.

We can see the two irreducible components of $X$ explicitly: define maps
\begin{align*}
\delta = q \times 1: S^{[n-1,n]}&\to S \times S^{[n-1,n]}&	\epsilon: S^{[n-2,n-1,n]}& \to S \times S^{[n-1,n]} \\
(\zeta, \xi) &\mapsto (\xi \setminus \zeta, \zeta, \xi)&	(\eta, \zeta, \xi)& \mapsto (\zeta \setminus \eta, \zeta, \xi).
\end{align*}
Then we have $\tilde\imath(X) = \im \delta \cup \im \epsilon$.  In \S\ref{irred} below we show that $S^{[n-2,n-1,n]}$, though not smooth \cite{cheah}, is indeed irreducible.  Let us manipulate the short exact sequence
\[ 0 \to \I_{\im \epsilon / \tilde\imath(X)} \to \O_{\tilde\imath(X)} \to \O_{\im \epsilon} \to 0. \]
The first term is isomorphic to the ideal sheaf of $\im \delta \cap \im \epsilon$ in $\im \delta$; moreover $\delta$ is an embedding, and $\im \delta \cap \im \epsilon = \delta(E)$, so the first term becomes $\delta_* \O(-E)$.  For the third term, note that the fiber of $\epsilon$ over a point $(x, \zeta, \xi) \in S \times S^{[n-1,n]}$ is a (possibly empty) projective space $\P \Hom(\O_x, \O_\zeta)^*$, so $\epsilon_* \O_{S^{[n-2,n-1,n]}} = \O_{\im \epsilon}$.  Thus we have
\[ 0 \to \delta_* \O(-E) \to \tilde\imath_* \O_X \to \epsilon_* \O_{S^{[n-2,n-1,n]}} \to 0. \]
Tensor with $\pi_2^* \O(E)$ and use the projection formula to get
\[ 0 \to \delta_* \O_{S^{[n-1,n]}} \to \tilde\imath_* \O_X \otimes \pi_2^* \O(E) \to \epsilon_* \epsilon^* \pi_2^* \O(E) \to 0. \]
Now apply $\pi_{SS*} \phi'_*$.  For the first term, observe that the diagram
\[ \xymatrix{
S^{[n-1,n]} \ar[d]_{\gamma = q \times f} \ar[rr]^-{\delta = q \times 1}
& & S \times S^{[n-1,n]} \ar[d]^{\phi' = 1 \times g \times q} \\
S \times S^{[n-1]} \ar[d]_{\pi_1}
& & S \times S^{[n]} \times S \ar[d]^{\pi_{SS}} \\
S \ar[rr]_-{\Delta} & & S \times S
} \]
commutes, and we have seen that $\gamma_* \O_{S^{[n-1,n]}} = \O_{S \times S^{[n-1]}}$, so the first term becomes $\O_\Delta \otimes \RGamma(\O_{S^{[n-1]}})$.  The second term becomes $R''F''$.  For the third term, observe that the composition
\[ S^{[n-2,n-1,n]} \xrightarrow\epsilon X \xrightarrow{\phi'} S \times S^{[n]} \times S \xrightarrow{\pi_{SS}} S \times S \]
sends a point $(\eta, \zeta, \xi)$ to $(\zeta\setminus\eta, \xi\setminus\zeta)$, hence equals the vertical composition in the diagram
\[ \xymatrix{
S^{[n-2,n-1,n]} \ar[d] \ar[r]^{\pi_2 \epsilon}
& S^{[n-1,n]} \ar[d]^\gamma \\
S^{[n-2,n-1]} \times S \ar[d] \ar[r]
& S^{[n-1]} \times S \\
S^{[n-2]} \times S \times S \ar[d] \\
S \times S.
} \]
Using base change around the square (again by Appendix \ref{base_change}) and the fact that $\gamma_* \O(E) = \O_{S^{[n-1]} \times S}$, we find that the third term becomes $\O_{S \times S} \otimes \RGamma(\O_{S^{[n-2]}})$.  Thus we get an exact triangle
\begin{equation} \label{for_q_eta_q_later}
\O_\Delta \otimes \RGamma(\O_{S^{[n-1]}}) \to R''F'' \to \O_{S \times S} \otimes \RGamma(\O_{S^{[n-2]}}),
\end{equation}
which must split because $\Ext^i(\O_{S \times S}, \O_\Delta) = H^i(\O_S)$ vanishes when $i$ is odd.

\subsection{Cancellation} \label{cancel}

Now we assemble what we know about $R'F'$, $R'F''$, $R''F'$, $R''F''$, and the maps between them to show that
\begin{align*}
RF &\cong \O_\Delta \otimes \RGamma(\O_{S^{[n-1]}}) \\
&= \O_\Delta \oplus \O_\Delta[-2] \oplus \dotsb \oplus \O_\Delta[-2n+2].
\end{align*}
For the reader's convenience we recall from the last two sections that
\begin{align*}
R'F' &= \O[2] \oplus \dotsb \oplus \O[-2n+2] \\
R'F'' &= \O[2] \oplus \dotsb \oplus \O[-2n+4] \\
R''F' &= \O \oplus \dotsb \oplus \O[-2n+2] \\
R''F'' &= \O \oplus \dotsb \oplus \O[-2n+4] \oplus \O_\Delta \oplus \dotsb \oplus \O_\Delta[-2n+2],
\end{align*}
where $\O$ is short for $\O_{S \times S}$.

We have a diagram of exact triangles
\[ \xymatrix{
R''F' \ar[r] \ar[d]
& R'F' \ar[r] \ar[d]
& RF' \ar[d] \\
R''F'' \ar[r]
& R'F'' \ar[r]
& RF''.
} \]
Let us take cohomology sheaves of this to get a diagram of exact sequences
\[ \xymatrix@C=1.5em{
\H^i(R''F') \ar[r] \ar[d]
& \H^i(R'F') \ar[r] \ar[d]
& \H^i(RF') \ar[r] \ar[d]
& \H^{i+1}(R''F') \ar[r] \ar[d]
& \H^{i+1}(R'F') \ar[d] \\
\H^i(R''F'') \ar[r]
& \H^i(R'F'') \ar[r]
& \H^i(RF'') \ar[r]
& \H^{i+1}(R''F'') \ar[r]
& \H^{i+1}(R'F'')
} \]
for various $i$.

\begin{itemize}

\item For $i=-2$, we have
\[ \xymatrix{
0 \ar[d] \ar[r]
& \O \ar@{=}[d] \ar[r]
& ? \ar[d] \ar[r]
& 0 \ar[d] \ar[r]
& 0 \ar[d] \\
0 \ar[r]
& \O \ar[r]
& ? \ar[r]
& 0 \ar[r]
& 0,
} \]
so
\[ \H^{-2}(RF') = \H^{-2}(RF'') = \O \]
and the map between them is an isomorphism.

\item For $i=-1$, we have
\[ \xymatrix{
0 \ar[r] \ar[d]
& 0 \ar[r] \ar[d]
& ? \ar[r] \ar[d]
& \O \ar@{=}[r] \ar[d]
& \O \ar@{=}[d] \\
0 \ar[r]
& 0 \ar[r]
& ? \ar[r]
& \O \oplus \O_\Delta \ar[r]
& \O.
} \]
Since the right-hand square is commutative, the map $\O \oplus \O_\Delta \to \O$ is split, so its kernel is $\O_\Delta$, so
\[ \H^{-1}(RF') = 0  \hspace{3em}  \H^{-1}(RF'') = \O_\Delta. \]

\item For $i=0$, we have
\[ \xymatrix{
\O \ar@{=}[r] \ar[d]
& \O \ar[r] \ar@{=}[d]
& ? \ar[r] \ar[d]
& 0 \ar[r] \ar[d]
& 0 \ar[d] \\
\O \oplus \O_\Delta \ar[r]
& \O \ar[r].
& ? \ar[r]
& 0 \ar[r]
& 0.
} \]
Then $\O \oplus \O_\Delta \to \O$ is surjective, so we get
\[ \H^0(RF') = \H^0(RF'') = 0. \]

\item For $1 \le i \le 2n-4$, we get the same result as for $i=-1$ and $i=0$ over and over.

\item For $i=2n-3$, we have
\[ \xymatrix{
0 \ar[r] \ar[d]
& 0 \ar[r] \ar[d]
& ? \ar[r] \ar[d]
& \O \ar@{=}[r] \ar[d]
& \O \ar[d] \\
0 \ar[r]
& 0 \ar[r]
& ? \ar[r]
& \O_\Delta \ar[r]
& 0,
} \]
so
\[ \H^{2n-3}(RF') = 0  \hspace{3em}  \H^{2n-3}(RF'') = \O_\Delta. \]

\item For $i=2n-2$, we have
\[ \xymatrix{
\O \ar@{=}[r] \ar[d]
& \O \ar[r] \ar[d]
& ? \ar[r] \ar[d]
& 0 \ar[r] \ar[d]
& 0 \ar[d] \\
\O_\Delta \ar[r]
& 0 \ar[r].
& ? \ar[r]
& 0 \ar[r]
& 0,
} \]
so
\[ \H^{2n-2}(RF') = \H^{2n-2}(RF'') = 0. \]
\end{itemize}

\noindent Now we take cohomology sheaves of the exact triangle
\[ RF \to RF' \to RF'' \]
to get a long exact sequence
\[ \begin{array}{rclcccccl}
0 & \to & \H^{-2}(RF) & \to & \O & = & \O \\
& \to & \H^{-1}(RF) & \to & 0 & \to & \O_\Delta \\
& \to & \H^0(RF) & \to & 0 & \dotsb & 0 \\
& \to & \H^{2n-3}(RF) & \to & 0 & \to & \O_\Delta \\
& \to & \H^{2n-2}(RF) & \to & 0 & \to & 0
\end{array} \]
which gives
\[ \H^i(RF) = \begin{cases}
\O_\Delta & i = 0, 2, \dotsc, 2n-2 \\
0 & \text{otherwise.}
\end{cases} \]
Thus $RF$ has a filtration whose associated graded object is $\O_\Delta \oplus \O_\Delta[-2] \oplus \dotsb \oplus \O_\Delta[-2n+2]$.  But $\Ext^i_{S \times S}(\O_\Delta, \O_\Delta) = HH^i(S)$ vanishes when $i$ is odd, so the filtration splits.

\subsection{Monad structure} \label{monad}

Having proved part (a) of Theorem \ref{hilb_calc_thm}, that
\begin{equation} \label{RF_iso}
RF \cong \O_\Delta \oplus \O_\Delta [-2] \oplus \dotsb \oplus \O_\Delta [-2n+2],
\end{equation}
we now consider the monad structure $RFRF \xrightarrow{R \epsilon F} RF$.\footnote{For background on monads in general see \cite[\S VI.1]{maclane}.}  Presumably it is like multiplication in $H^*(\P^{n-1})$, but we will prove the following weaker statement, which is sufficient for our purposes in \S\ref{P-functors}:
{\renewcommand \thethm {\ref{hilb_calc_thm}}
\addtocounter{thm}{-1}
\begin{thm}
(b) The isomorphism \eqref{RF_iso} can be chosen so that the map
\[ RF[-2] \hookrightarrow RFRF  \xrightarrow{R \epsilon F} RF, \]
when written in components
\begin{multline*}
\O_\Delta[-2] \oplus \O_\Delta[-4] \oplus \dotsb \oplus \O_\Delta[-2n] \\
\longrightarrow \O_\Delta \oplus \O_\Delta[-2] \oplus \dotsb \oplus \O_\Delta[-2n+2],
\end{multline*}
is of the form
\[ \begin{pmatrix}
0  \\
1 & 0 \\
  & 1 & 0 \\
& & \ddots & \ddots \\
& & & 1 & 0 & * \\
& & & & 1 & *
\end{pmatrix}. \]
\end{thm}}

We introduce the endofunctor
\[ \Phi := \id_S \otimes \RGamma(\O_{S^{[n]}}) \]
of $D^b(S)$, with a monad structure given by the ring structure in the second factor, and a map of monads $\varphi\colon \Phi \to RF$.  We define $\varphi$ as the adjoint to the natural transformation $\tilde\varphi\colon F \otimes \RGamma(\O_{S^{[n]}}) \to F$ given by
\[ F(-) \otimes \O_{S^{[n]}} \otimes \RGamma(\O_{S^{[n]}}) \xrightarrow{1 \otimes \ev} F(-) \otimes \O_{S^{[n]}}. \]
That is, $\varphi$ is the composition
\[ \id_S \otimes \RGamma(\O_{S^{[n]}})
\xrightarrow{\eta \otimes 1}
RF \otimes \RGamma(\O_{S^{[n]}})
\xrightarrow{R \tilde\varphi}
RF. \]
It is straightforward to check that this is a map of monads.

\begin{lem} \label{monad_claim_1}
The map $\varphi\colon\Phi \to RF$ induces an isomorphism on $\H^i$ for $i = 0, 2, \dotsc, 2n-2$.
\end{lem}

Before proving this claim we show how it implies Theorem \ref{hilb_calc_thm}(b).  Write $\varphi$ in components
\begin{multline*}
\O_\Delta \oplus \O_\Delta[-2] \oplus \dotsb \oplus \O_\Delta[-2n+2] \oplus \O_\Delta[-2n] \\
\longrightarrow \O_\Delta \oplus \O_\Delta[-2] \oplus \dotsb \oplus \O_\Delta[-2n+2],
\end{multline*}
so it is of the form
\begin{equation} \label{phi_a_priori}
\varphi = \begin{pmatrix}
a_0 & * & * \\
& a_1 & * & * \\
& & \ddots & \ddots & \ddots \\
& & & a_{n-2} & * & * \\
& & & & a_{n-1} & *
\end{pmatrix}
\end{equation}
where $a_0, a_1, \dotsc, a_{n-1} \in \Hom(\O_\Delta, \O_\Delta) = \C$ are non-zero.  Then we can compose the isomorphism \eqref{RF_iso} with an automorphism of 
\[ \O_\Delta \oplus \O_\Delta [-2] \oplus \dotsb \oplus \O_\Delta [-2n+2] \]
of the form
\[ \begin{pmatrix}
a_0^{-1} & * & \cdots & * & * \\
& a_1^{-1} & \cdots & * & * \\
& & \ddots & \vdots & \vdots \\
& & & a_{n-2}^{-1} & * \\
& & & & a_{n-1}^{-1}
\end{pmatrix} \]
so that \eqref{phi_a_priori} becomes
\begin{equation} \label{phi_normalized}
\varphi = \begin{pmatrix}
1 \\
& 1 \\
& & \ddots \\
& & & 1 & & * \\
& & & & 1 & *
\end{pmatrix}.
\end{equation}
Now we have a commutative diagram
\begin{equation} \label{Phi_dominates_RF}
\begin{split} \xymatrix{
\Phi[-2] \ar@{^{(}->}[r] \ar[d]_{\varphi[-2]}
& \Phi\Phi \ar[r] \ar[d]_{\varphi\varphi}
& \Phi \ar[d]_\varphi \\
RF[-2] \ar@{^{(}->}[r]
& RFRF \ar[r]
& RF;
} \end{split} \end{equation}
to see that the left-hand square commutes, observe that
\[ \xymatrix{
\O_\Delta[-2] \ar@{^{(}->}[r] \ar@{=}[d]
& \Phi \ar[d]_\varphi \\
\O_\Delta[-2] \ar@{^{(}->}[r]
& RF
} \]
commutes because $\varphi$ is of the form \eqref{phi_normalized}.  Now the composition across the top of \eqref{Phi_dominates_RF} is
\[ \begin{pmatrix}
0 \\
1 & 0 \\
& 1 & 0 \\
& & \ddots & \ddots \\
& & & 1 & 0 \\
\end{pmatrix}, \]
and the outside vertical maps are \eqref{phi_normalized}, so the composition across the bottom is necessarily
\[ \begin{pmatrix}
0  \\
1 & 0 \\
& \ddots & \ddots\\
& & 1 & 0 & * \\
& & & 1 & *
\end{pmatrix}, \]
as desired.

\bigskip
Now we work toward proving Lemma \ref{monad_claim_1}.  We calculated $RF$ by first calculating $R''F''$ and then chasing through some long exact sequences, so to understand the map $\varphi\colon \Phi \to RF$ we will first study the analogous map $\varphi''\colon \Phi \to R''F''$.  Recall that
\[ \H^i(R''F'') = \begin{cases}
\O \oplus \O_\Delta & i = 0, 2, \dotsc, 2n-4 \\
\O_\Delta & i = 2n-2 \\
0 & \text{otherwise.}
\end{cases} \]

\begin{lem} \label{monad_claim_2}
The map $\varphi''\colon \Phi \to R''F''$ induces a non-zero map on $\H^i$ for $i = 0, 2, \dotsc, 2n-2$.
\end{lem}

\begin{proof}
Again let $q$ and $g$ be as in
\[ \xymatrix{
S^{[n-1,n]} \ar[d]_q \ar[r]^g & S^{[n]} \\
S
} \]
and recall that $F'' = g_* q^*$, so $R'' = q_* g^!$.  We can factor $\varphi''$ as
\begin{equation} \label{factorization_of_phi''}
\Phi \xrightarrow\psi q_* q^* \xrightarrow{q_* \eta q^*} q_* g^! g_* q^*,
\end{equation}
as follows.  Let $\psi$ be adjoint to the map $\tilde\psi\colon q^* \otimes R\Gamma(\O_{S^{[n]}}) \to q^*$ given by
\begin{equation} \label{psi_tilde}
q^*(-) \otimes g^* \O_{S^{[n]}} \otimes R\Gamma(\O_{S^{[n]}}) \xrightarrow{1 \otimes g^*{\ev}} q^* (-) \otimes g^* \O_{S^{[n]}}.
\end{equation}
More explicitly, for $i \in \Z$ and $\tau \in H^i(\O_{S^{[n]}}) = \Hom(\O_{S^{[n]}}[-i], \O_{S^{[n]}})$ we are talking about
\[ q^*(-) \otimes g^* \O_{S^{[n]}}[-i] \xrightarrow{1 \otimes g^* \tau} q^*(-) \otimes g^* \O_{S^{[n]}}. \]
Apply $g_*$ to \eqref{psi_tilde} and use functoriality of the projection formula to get
\[ g_* q^*(-) \otimes \O_{S^{[n]}} \otimes R\Gamma(\O_{S^{[n]}})  \xrightarrow{1 \otimes \ev} g_* q^*(-) \otimes \O_{S^{[n]}}, \]
which is exactly the map
\[ \tilde\varphi''\colon F'' \otimes \RGamma(\O_{S^{[n]}}) \to F'' \]
adjoint to $\varphi''$.  So $g_* \tilde\psi = \tilde\varphi''$, and thus the 
diagram
\[ \xymatrix{
\id_S \otimes \RGamma(\O_{S^{[n]}})
\ar[r]^-{\eta \otimes 1}
\ar[rd]_{\eta \otimes 1}
& q_* q^* \otimes \RGamma(\O_{S^{[n]}})
\ar[r]^-{q_* \tilde\psi}
\ar[d]^{q_* \eta q^* \otimes 1}
& q_* q^*
\ar[d]^{q_* \eta q^*} \\
& q_* g^! g_* q^* \otimes \RGamma(\O_{S^{[n]}})
\ar[r]^-{q_* g^! \tilde\varphi''}
& q_* g^! g_* q^*
} \]
commutes, which gives the desired factorization \eqref{factorization_of_phi''}.
\bigskip

Recall that
\begin{multline*}
q_* q^* = \id_S \otimes q_* \O_{S^{[n-1,n]}} = \id_S \otimes \RGamma(\O_{S^[n-1]}) \\
= \O_\Delta \oplus \O_\Delta[-2] \oplus \dotsb \oplus \O_\Delta[-2n+2].
\end{multline*}
From \eqref{for_q_eta_q_later} we know that the cone on $q_* q^* \xrightarrow{q_* \eta q^*} R''F''$ is $\O \oplus \O[-2] \oplus \dotsb \oplus \O[-2n+4]$, so the induced map $\H^i(q_* q^*) \to \H^i(R''F'')$ is certainly non-zero for $i = 0, 2, \dotsc, 2n-2$.  Thus we need only show that $\psi\colon \Phi \to q_* q^*$ induces isomorphisms on $\H^i$ for $i = 0, 2, \dotsc, 2n-2$.  But this boils down to the claim that the natural map
\[ H^i(\O_{S^{[n]}}) \to H^i(q_* g^* \O_{S^{[n]}}) \]
is an isomorphism for $i = 0, 2, \dotsc, 2n-2$, which we proved in \S\ref{easy}.
\end{proof}

\begin{proof}[Proof of Lemma \ref{monad_claim_1}]
First we prove the claim for $i = 0, 2, \dotsc, 2n-4$, postponing the case of $i=2n-2$ for a moment.  The diagram
\[ \xymatrix{
\Phi \ar[d]_\varphi \ar[r]^{\varphi''} & R''F'' \ar[d] \\
RF \ar[r] & R''F[1]
} \]
commutes for formal reasons.  To understand $R''F[1]$, take the triangle
\[ R'F \to RF \to R''F[1] \]
and recall that $R'F = \O[-2n+2]$.  Thus the map $\H^i(RF) \to \H^i(R''F[1])$ is an isomorphism for $i = 0,2,\dotsc,2n-4$, so it is enough to show that $\H^i(\Phi) \to \H^i(R''F[1])$ is an isomorphism for $i = 0, 2, \dotsc, 2n-4$.  From the diagram chase in \S\ref{cancel} we see that $\H^i(R''F'') \to \H^i(R''F[1])$ for $i = 0, 2, \dotsc, 2n-4$ is the projection $\O \oplus \O_\Delta \to \O_\Delta$, which together with Lemma \ref{monad_claim_2} gives what we want.

Now we prove the claim for $i = 2n-2$.  It is not obvious whether the map $\H^{2n-2}(R'F) \to \H^{2n-2}(RF)$, which is $\O \to \O_\Delta$, is zero or non-zero, but rather than decide the question we give a proof in both cases.

\smallskip
Case 1: $\H^{2n-2}(R'F) \to \H^{2n-2}(RF)$ is zero.  Then we see that the map $\H^{2n-2}(RF) \to \H^{2n-2}(R''F[1])$ is an isomorphism.  Since $\H^{2n-1}(R''F') = 0$ we see that the map $\H^{2n-2}(R''F'') \to \H^{2n-2}(R''F[1])$, which is $\O_\Delta \to \O_\Delta$, is an isomorphism as well, which together with Lemma \ref{monad_claim_2} gives what we want.

\smallskip
Case 2: $\H^{2n-2}(R'F) \to \H^{2n-2}(RF)$ is non-zero.  Then we have
\[ \H^i(R''F[1]) = \begin{cases}
\O_\Delta & i = 0, 2, \dotsc, 2n-4 \\
I_\Delta & i = 2n-3 \\
0 & \text{otherwise.}
\end{cases} \]
Since $\Ext^i(\O_\Delta, \O_\Delta)$ vanishes when $i$ is even and $\Ext^i(I_\Delta, \O_\Delta)$ vanishes when $i$ is odd, $R''F[1]$ must split as the sum of its cohomology sheaves:
\[ R''F[1] \cong \O_\Delta \oplus \O_\Delta[-2] \oplus \dotsb \oplus \O_\Delta[-2n+4] \oplus I_\Delta[-2n+3]. \]
Moreover the component $\O_\Delta[-2n+2] \to I_\Delta[-2n+3]$ of $RF \to R''F[1]$ is non-zero.

Observe that $\Hom(\O_\Delta[-2n+2], I_\Delta[-2n+3]) = \Ext^1(\O_\Delta, I_\Delta) = \C$.  Thus it is enough to show that the component $\O_\Delta[-2n+2] \to I_\Delta[-2n+3]$ of $R''F'' \to R''F[1]$ is non-zero.  Consider the triangle
\[ R''F' \to R''F'' \to R''F[1] \]
and recall that
\[ R''F' = \O \oplus \O[-2] \oplus \dotsb \oplus \O[-2n+2]. \]
Since $\H^{2n-2}(R''F[1]) = 0$, the component $\O[-2n+2] \to \O_\Delta[-2n+2]$ of $R''F' \to R''F''$ is non-zero, so the component $\O_\Delta[-2n+2] \to I_\Delta[-2n+3]$ of $R''F'' \to R''F[1]$ is non-zero as desired.
\end{proof}

\let \oldthesubsection \thesubsection
\renewcommand \thesubsection {\arabic{section}.A}
\subsection{Appendix: Irreducibility of \texorpdfstring{$S^{[n-2,n-1,n]}$}{S[n-2,n-1,n]}}\label{irred}
\let \thesubsection \oldthesubsection

In this section we show that if $S$ is a smooth surface then the nested Hilbert scheme $S^{[n-2,n-1,n]}$ is irreducible, which we needed in \S\ref{hard}.  Note that not all nested Hilbert schemes are irreducible: if $n \gg 0$ then $S^{[1,2,\dotsc,n]}$ has components whose dimension is greater than the expected $2n$.\footnote{I thank Mark Haiman for explaining this to me.}

To lighten the notation we work with $S^{[n-1,n,n+1]}$.  Recall that $S^{[n]}$ and $S^{[n-1,n]}$ are smooth of dimension $2n$ for all $n$.  Write $S^{[n-1,n,n+1]}$ as the intersection
\[ (S^{[n-1,n]} \times S^{[n+1]}) \cap (S^{[n-1]} \times S^{[n,n+1]}) \subset S^{[n-1]} \times S^{[n]} \times S^{[n+1]}. \]
Because the ambient space is smooth, every component of the intersection has at least the expected dimension
\[ \big(2n + (2n+2)\big) + \big((2n-2) + (2n+2)\big) - \big((2n-2) + 2n + (2n+2)\big) = 2n+2. \]

Consider the fiber square
\[ \xymatrix{
S^{[n-1,n,n+1]} \ar[r]^-{\tilde g} \ar[d] & S^{[n,n+1]} \ar[d]^-f \\
S^{[n-1,n]} \ar[r]_-g & S^{[n]}.
} \]
Let $U = S^{[n,n+1]} \setminus E$, where $E$ is the exceptional divisor.  Then the fibers of $f|_U$ are the fibers of $(S^{[n]} \times S) \setminus Z \to S^{[n]}$, which are irreducible surfaces ($S$ minus finitely many points), so $\tilde g^{-1}(U) = S^{[n-1,n]} \times_{S^{[n]}} U$ is irreducible.  Thus it is enough to show that $\dim \tilde g^{-1}(E) < 2n+2$.

The maps $g\colon S^{[n-1,n]} \to S^{[n]}$ and $f|_E\colon E \to S^{[n]}$ factor through $Z$.  Ellingsrud and Str{\o}mme \cite[\S3]{es} partition $Z$ into locally closed subsets
\[ W_i = \{ (x,\zeta) \in Z : \dim (\I_\zeta|_x) = i \} \hspace{3em} i \ge 2 \]
and show that the fiber of $E \to Z$ over $W_i$ is $\P^{i-1}$, the fiber of $S^{[n-1,n]} \to Z$ over $W_i$ is $\P^{i-2}$, and $\dim W_i \le 2n + 4 - 2i$.  Now $\tilde g^{-1}(E) = S^{[n-1,n]} \times_{S^{[n]}} E$ maps to $Z \times_{S^{[n]}} Z$, which we partition into
\[ W_i \times_{S^{[n]}} W_j = \{ (w,w') \in W_i \times W_j : \pi(w) = \pi(w') \}, \]
where $\pi: Z \to S^{[n]}$.  Since $Z \times_{S^{[n]}} Z \to S^{[n]}$ is finite, the dimension of $W_i \times_{S^{[n]}} W_j$ is the same as that of its image $\pi(W_i) \cap \pi(W_j)$, which is at most $2n+4-2\max\{i,j\}$.  Thus the preimage of $W_i \times_{S^{[n]}} W_j$ in $S^{[n-1,n]} \times_{S^{[n]}} E$ has dimension at most
\[ (2n+4-2\max\{i,j\})+(i-2)+(j-1) \le 2n+1, \]
which gives the desired result.

%% file: P-functors.tex

\section{\texorpdfstring{$\P$}{P}-Functors} \label{P-functors}

\subsection{Definition}

In view of Theorem \ref{hilb_calc_thm}, we need to define $\P$-functors, generalizing Huybrechts and Thomas's $\P$-objects \cite{ht}.  Let $X$ be a smooth $2n$-dimensional complex projective variety, and recall that an object $\E \in D^b(X)$ is called a \emph{$\P^n$-object} if $\Ext^*(\E, \E) \cong H^*(\P^n, \C)$ as rings, and $\E \otimes \omega_X \cong \E$.  The first example is a line bundle on a hyperk\"ahler variety.  The second is the structure sheaf of a Lagrangian $\P^n$ in a hyperk\"ahler variety: for example, $\P^n$ sitting in the total space of its cotangent bundle, or if $S$ is a K3 surface containing a rational curve $C \cong \P^1$ then $C^{[n]} \cong \P^n \subset S^{[n]}$.

\begin{defn}
A \emph{$\P^n$-functor} is a functor $F\colon \A \to \B$ with adjoints $L$ and $R$ such that:
\begin{enumerate}
\renewcommand \theenumi {\alph{enumi}}
\item There is an autoequivalence $H$ of $\A$ such that
\begin{equation} \label{ring_iso} 
RF \cong \id \oplus H \oplus H^2 \oplus \dotsb \oplus H^n.
\end{equation}
\item \label{ring_str} The map
\[ HRF \hookrightarrow RFRF \xrightarrow{R \epsilon F} RF, \]
when written in components
\[ H \oplus H^2 \oplus \dotsb \oplus H^n \oplus H^{n+1} \to \id \oplus H \oplus H^2 \oplus \dotsb \oplus H^n, \]
is of the form
\[ \begin{pmatrix}
* & * & \cdots & * & * \\
1 & * & \cdots & * & * \\
0 & 1 & \cdots & * & * \\
\vdots & \vdots & \ddots & \vdots & \vdots \\
0 & 0 & \cdots & 1 & * \\
\end{pmatrix}. \]
(This models the fact that the map $\C[h]/h^{n+1} \to \C[h]/h^{n+1}$ given by multiplication by $h$ has kernel $\C \cdot h^n$ and cokernel $\C \cdot 1$.)
\item $R \cong H^n L$.  If $\A$ and $\B$ have Serre functors, this is equivalent to $S_\B F H^n \cong F S_\A$.
\end{enumerate}
As in \S\ref{spherical} we need all functors to be induced by Fourier--Mukai kernels, and the isomorphism \eqref{ring_iso} to be induced by a map of kernels.
\end{defn}

Cautis \cite{cautis_about} has recently made a similar definition; he considers only $H = [-2]$, but we will see interesting examples with $H = [-1]$ as well.  Of course one would like to find examples with more exciting $H$.

After giving examples of $\P$-functors, we construct the $\P$-twist associated to a $\P$-functor.  There is some work to do beyond simply quoting Huybrechts and Thomas, because in taking a certain double cone we face a choice that they do not.

\subsection{Examples}

\begin{enumerate}
\item If $\E \in D^b(X)$ is a $\P^n$-object then $F = \E \otimes -\colon D^b(\text{point}) \to D^b(X)$ is a $\P^n$-functor with $H = [-2]$.

\item The functor $F\colon D^b(S) \to D^b(S^{[n]})$ in the previous section is a $\P^{n-1}$-functor with $H = [-2]$.  The condition $S_{S^{[n]}} F H^{n-1} \cong F S_S$ is satisfied because $S_{S^{[n]}} = [2n]$ and $S_S = [2]$.

\item \label{split_spherical} A split spherical functor $F\colon \A \to \B$, that is, one where the exact triangle
\[ \id_\A \xrightarrow{\eta} RF \to C \]
is split, so $RF \cong \id_\A \oplus C$, is a $\P^1$-functor with $H = C$.  The $\P^1$-twist that we will construct in \S\ref{def_of_P} will coincide with the square of the spherical twist, just as in \cite[Prop.~2.9]{ht}.

\item Let $q\colon E \to Z$ be a $\P^n$-bundle, and let $i\colon E \to \Omega^1_q$ be the zero section of the relative cotangent bundle.  Then $F := i_* q^*$ is a $\P^n$-functor with $H = [-2]$, as follows.  The normal bundle of $E$ is $\Omega^1_q$, so $\omega_i = \omega_q$, so
\begin{align*}
R &= q_* i^! \\
&= q_*(\omega_i[-n] \otimes i^* -) \\
&= q_*(\omega_q[n][-2n] \otimes i^* -) \\
&= q_! i^*[-2n] \\
&= H^n L.
\end{align*}
Let $p\colon \Omega^1_q \to E$ be the projection; then for any $\F \in D^b(E)$ we have
\begin{align*}
i^! i_* \F &= i^! i_* i^* p^* \F \\
&= i^!(p^* \F \otimes i_* \O_E) \\
&= i^* p^* \F \otimes i^! i_* \O_E \\
&= \F \otimes i^! i_* \O_E.
\end{align*}
Moreover,
\begin{align*}
i^! i_* \O_E &= \omega_i[-n] \otimes i^* i_* \O_E \\
&= \omega_q[-n] \otimes (\O_E \oplus T_q[1] \oplus \Lambda^2 T_q[2] \oplus \dotsb \oplus \Lambda^{n} T_q[n]) \\
&= \O_E \oplus \Omega^1_q[-1] \oplus \Omega^2_q[-2] \oplus \dotsb \oplus \Omega^n_q[-n].
\end{align*}
Thus for any $\G \in D^b(Z)$ we have
\begin{align*}
RF \G &= q_* i^! i_* q^* \G \\
&= q_*(q^* \G \otimes i^! i_* \O_E) \\
&= \G \otimes q_* (\O_E \oplus \Omega^1_q[-1] \oplus \Omega^2_q[-2] \oplus \dotsb \oplus \Omega^n_q[-n]) \\
&= \G \otimes (\O_Z \oplus \O_Z[-2] \oplus \O_Z[-4] \oplus \dotsb \oplus \O_Z[-2n])
\end{align*}
since $q_* \Omega^k_q = \O_Z[-k]$.

If $Z$ is a point then the isomorphism $\Ext^*(\O_E, \O_E) \cong H^*(\P^n)$ is well-known to be a ring isomorphism.  For more general $Z$, note that the zeros below the diagonal in hypothesis (\ref{ring_str}) come for free, since $\Ext^{<0}(\O_\Delta, \O_\Delta) = 0$.  Thus we need only check that $RF[-2] \to RFRF \xrightarrow{R\epsilon F} RF$ induces an isomorphisms on $\H^i$ for $2 \le i \le 2n-2$.  This can be done pointwise, where it follows from the well-known case.

\item The following example is due to Kawamata \cite{kawamata}.  Let $X$ be a 3-fold with an $A_2$ singularity, $\tilde X \to X$ the blowup of the singular point, $E$ the exceptional divisor, $l$ a ruling of the quadric cone $E$, and $\B = {}^\perp\langle \mathcal{O}_E(E) \rangle \subset  D^b(\tilde X)$.  Then he shows that $\operatorname{Perf}(X)^\perp \subset \B$ is generated by one object $\E = \O_E(-l)$, that $\Ext^*(\E, \E) \cong \C[h]/h^3$ as graded rings with $\deg h = 1$, and that the Serre functor $S_\B$ acts as $S_\B \E \cong \E[2]$.  Thus the functor $\E \otimes -\colon D^b(\text{point}) \to \B$ is a $\P^2$-functor with $H = [-1]$.

\item \label{toda_kawamata} Examples like the previous one, which we might sheepishly call $\mathbb{RP}^n$-objects, are equivalent to Toda's fat spherical objects with $A = \C[\epsilon]/\epsilon^2$ the ring of dual numbers, as follows.  Let $\E \in D^b(X)$ be an object such that $\Ext^*(\E,\E) = \C[h]/h^{n+1}$ with $\deg h = 1$ instead of 2.  Let $\E' \in D^b(\Spec A \times X)$ be the first-order deformation corresponding to $h \in \Ext^1(\E,\E)$, and let $\pi\colon \Spec A \times X \to X$; then $\pi_* \E'$ is the non-trivial extension
\[ 0 \to \E \to \pi_* \E' \to \E \to 0. \]
Apply $\Hom(-, \E)$ to get
\begin{align*}
0 &\to \Hom(\E,\E) \to \Hom(\pi_* \E', \E) \to \Hom(\E,\E) \\
&\xrightarrow{\cdot h} \Ext^1(\E,\E) \to \Ext^1(\pi_* \E', \E) \to \Ext^1(\E,\E) \\
&\xrightarrow{\cdot h} \dotsb \\
&\xrightarrow{\cdot h} \Ext^n(\E,\E) \to \Ext^n(\pi_* \E', \E) \to \Ext^n(\E,\E) \to 0.
\end{align*}
Then the boundary maps are all isomorphisms, so $\Ext^*(\pi_* \E', \E) \cong H^*(S^n, \C)$, so $\E'$ is a fat spherical object.  From the same long exact sequence we see that the converse holds as well: if $\Ext^*(\pi_* \E', \E) \cong H^*(S^n, \C)$ then $\Ext^*(\E,\E) \cong \C[h]/h^{n+1}$ as rings.  The $\P^n$-twist that we will construct in \S\ref{def_of_P} will coincide with the fat spherical twist associated to $\E^*$.

Thus, for example, let $X$ be a 3-fold and let $C \subset X$ be a $(0,-2)$-curve which deforms to first order but not to second order; then the functor $\O_C \otimes -\colon D^b(\text{point}) \to D^b(X)$ is a $\P^3$-functor with $H = [-1]$.

\end{enumerate}

\subsection{Construction of the \texorpdfstring{$\P$}{P}-twist} \label{def_of_P}

We first recall Huybrechts and Thomas's definition of the $\P$-twist associated to a $\P$-object $\E$.  Let $h\colon \E[-2] \to \E$ be the map corresponding to a generator of $\Ext^2(\E, \E)$, and let $h^*\colon \E^*[-2] \to \E^*$ be its transpose.  Then $P\colon D^b(X) \to D^b(X)$ is the functor induced by the double cone
\[ \cone(\cone(\E^* \boxtimes \E[-2] \xrightarrow{h^* \boxtimes \id - \id \boxtimes h} \E^* \boxtimes \E) \xrightarrow{\operatorname{tr}} \O_\Delta) \]
in $D^b(X \times X)$.  Since $\Ext^{-1}(\E^* \boxtimes \E[-2], \O_\Delta) = \Ext^1(\E, \E) = 0$, there is a unique way to take this double cone.
\bigskip

To define the $\P$-twist associated to a $\P$-functor $F$, first let $j\colon H \to RF$ be the map coming from the splitting \eqref{ring_iso}, and let $f$ be the composition
\[ FHR \xrightarrow{F j R} FRFR \xrightarrow{\epsilon FR - FR\epsilon} FR, \]
where we recall that $\epsilon\colon FR \to \id_\B$ is the counit of the adjunction.  Then $f$ replaces $h^* \boxtimes \id - \id \boxtimes h$.  The composition
\[ FHR \xrightarrow{f} FR \xrightarrow\epsilon \id_\B \]
is zero, so we can take the double cone
\begin{equation} \label{double_cone}
\cone(\cone(FHR \xrightarrow{f} FR) \xrightarrow\epsilon \id_\B).
\end{equation}
But there need not be a unique way to take this double cone, since
\[ \Ext^{-1}(FHR, \id_\B) = \Ext^{-1}(H, RF) \]
need not vanish: for example, if $H = [-2]$ then this $\Ext$ group is $HH^1(\A)$, and if $H = [-1]$ it is $HH^0(\A)$, which never vanishes.

We will make an explicit choice for the double cone.  The functors $T = \cone \epsilon$ and $C = \cone \eta$, which were equivalences when $F$ was spherical, will now be used in an auxiliary way.  We will produce a lift $\tilde f$ as in the diagram
\begin{equation} \label{f_tilde}
\begin{split}
\xymatrix{
& T[-1] \ar[d] \\
FHR \ar@{..>}[ru]^{\tilde f} \ar[r]_f & FR \ar[d]^\epsilon \\
& \id_\B 
}
\end{split}
\end{equation}
and define $P = \cone \tilde f[1]$; using the octahedral axiom one can check that this is the same as \eqref{double_cone}.

To produce $\tilde f$, we will lift
\begin{equation} \label{lift}
\xymatrix{
& & T[-1] \ar[d] \\
FRFR \ar[rr]_{\epsilon FR - FR\epsilon} \ar@{..>}[rru] & & FR
}
\end{equation}
using the splitting of $FRF$ discussed in \S\ref{FRF}.  Consider the diagram
\[ \xymatrix{
TFR[-1] \ar@/_{2em}/[dd]_[@]\cong \ar[d] \ar[rr]^{-T\epsilon[-1]} & & T[-1] \ar[d] \\
FRFR \ar[d] \ar[rr]_{\epsilon FR - FR \epsilon} & & FR \\
FCR.
} \]
The identity $FRFR \to FRFR$ is the sum of the two idempotents
\begin{gather}
FRFR \to FCR \xleftarrow\cong TFR[-1] \to FRFR \label{first_idempotent} \\
FRFR \xrightarrow{\epsilon FR} FR \xrightarrow{F\eta R} FRFR. \label{second_idempotent}
\end{gather}
The composition of \eqref{second_idempotent} with
\begin{equation} \label{eFR_FRe}
FRFR \xrightarrow{\epsilon FR - FR\epsilon} FR
\end{equation}
is zero, so \eqref{eFR_FRe} equals the composition of \eqref{first_idempotent} with \eqref{eFR_FRe}, that is,
\[ FRFR \to FCR \xleftarrow\cong TFR[-1] \to FRFR \xrightarrow{\epsilon FR - FR\epsilon} FR. \]
Considering the diagram again we see that this equals
\[ FRFR \to FCR \xleftarrow\cong TFR[-1] \xrightarrow{-T\epsilon[-1]} T[-1] \to FR, \]
so for the lift in \eqref{lift} we can take
\[ FRFR \to FCR \xleftarrow\cong TFR[-1] \xrightarrow{-T\epsilon[-1]} T[-1]. \]

\begin{defn}
If $F$ is a $\P$-functor, the associated \emph{$\P$-twist} is the cone $P$ on the following composition:
\[ P := \cone(FHR[1] \xrightarrow{FjR[1]} FRFR[1] \to FCR[1] \xleftarrow\cong TFR \xrightarrow{-T\epsilon} T). \]
\end{defn}

In Example \ref{split_spherical} above we claimed that for $\P^1$-functors, which are the same as split spherical functors, the $\P^1$-twist is the square of the spherical twist.  To see this, observe that the composition $FHR \to FRFR \to FCR$ is an isomorphism in this case, so
\[ P = \cone(TFR \xrightarrow{-T\epsilon} T) = TT \]
since $T = \cone \epsilon$.
\bigskip

In Example \ref{toda_kawamata} we claimed that for an ``$\mathbb{RP}^n$-object'' $\E$, the fat spherical twist associated to $\E$ coincides with the $\P^n$-twist associated to $\E^*$.  Let $\E$ and $\E'$ be as in that example and let $h\colon \E \to \E[1]$ correspond to $h \in \Ext^1(\E,\E)$, so we have an exact triangle
\[ \E \to \pi_* \E' \to \E \xrightarrow{h} \E[1]. \]
Let $h^*\colon \E[-1] \to \E$ be its transpose, so we have an exact triangle
\[ \E^* \to \pi_*(\E'^*) \to \E^* \xrightarrow{-h^*} \E^*[1], \]
where the minus sign is due to the shift.  If $F\colon D^b(\Spec A) \to D^b(X)$ is the functor induced by $\E'$ and $R$ is its right adjoint, then $FR$ is induced by
\[ (\pi_{13*}(\pi_{12}^* \E'^* \otimes \pi_{23}^* \E'))^* \in D^b(X \times X), \]
where $\pi_{ij}$ are the projections from $X \times \Spec A \times X$.  Now $\pi_{23}^* \E' \in D^b(X \times \Spec A \times X)$ is a first-order deformation of $\O_X \boxtimes \E$, and $\pi_{12}^* \E'^*$ is a deformation of $\E^* \boxtimes \O_X$, so $\pi_{12}^* \E'^* \otimes \pi_{23}^* \E'$ is a deformation of $\E^* \boxtimes \E$, and when we push down we get an exact triangle
\[ \E^* \boxtimes \E \to \underbrace{\pi_{13*}(\pi_{12}^* \E'^* \otimes \pi_{23}^* \E')}_{= (FR)^*} \to \E^* \boxtimes \E \xrightarrow{\id \boxtimes h - h^* \boxtimes \id} \E^* \boxtimes \E[1]. \]
Thus we have
\begin{align*}
T &= \cone(FR \to \id) \\
&= \cone(\cone(\E \boxtimes \E^*[-1] \xrightarrow{\id \boxtimes h^* - h \boxtimes \id} \E \boxtimes \E^*) \to \id)
\end{align*}
as claimed.

\subsection{Action on the spanning class and cohomology}

If $\E$ is a $\P^n$-object then Huybrechts and Thomas show that the $\P^n$-twist $P$ sends $\E$ to $\E[-2n]$ and acts as the identity on $\E^\perp$.  The action on cohomology is trivial: the cone on
\[ \E^* \boxtimes \E[-2] \to \E^* \boxtimes \E \]
is zero in K-theory, so $\O_\Delta$ and $P$ are the same in K-theory.

We generalize this $\P^n$-functors.  Again the $\P$-twist
\[ P = \cone(FHR[1] \to T) \]
acts as the identity on $(\im F)^\perp = \ker R$: if $RB = 0$ then $PB = TB = B$.  It acts on $\im F$ as follows:
\begin{prop} $PF \cong FH^{n+1}[2]$. \end{prop}
\begin{proof}
The functor $PF[-1]$ is the cone on
\[ FHRF \xrightarrow{\tilde f F} TF[-1], \]
where $\tilde f$ is as in \eqref{f_tilde}.  If we post-compose with the isomorphism
\[ TF[-1] \to FRF \to FC \]
then the cone is unchanged, so $PF[-1]$ is the cone on
\[ FHRF \xrightarrow{fF} FRF \to FC, \]
or in more detail
\begin{equation} \label{cone_to_take}
FHRF \xrightarrow{FjRF} FRFRF \xrightarrow{\epsilon FRF - FR\epsilon F} FRF \to FC.
\end{equation}
Let us write \eqref{cone_to_take} in components with respect to the decompositions
\begin{align*}
FHRF &= FH \oplus FH^2 \oplus \dotsb \oplus FH^n \oplus FH^{n+1} \\
FC &= FH \oplus FH^2 \oplus \dotsb \oplus FH^n.
\end{align*}
By hypothesis (\ref{ring_str}) of the definition, the composition
\[ FHRF \xrightarrow{FjRF} FRFRF \xrightarrow{FR\epsilon F} FRF \to FC \]
is of the form
\[ \begin{pmatrix}
1 & * & \cdots & * & * \\
0 & 1 & \cdots & * & * \\
\vdots & \vdots & \ddots & \vdots & \vdots \\
0 & 0 & \cdots & 1 & * \\
\end{pmatrix}. \]
On the other hand, we can get the composition
\[ FHRF \xrightarrow{FjRF} FRFRF \xrightarrow{\epsilon FRF} FRF \to FC \]
by taking
\[ FH \xrightarrow{Fj} FRF \xrightarrow{\epsilon F} F, \]
applying $RF = \id_B \oplus H \oplus \dotsb \oplus H^n$ on the right, and post-composing with $FRF \to FC$; hence it is of the form
\[ \begin{pmatrix}
0 & * \\
& 0 & * \\
& & \ddots & \ddots \\
& & & 0 & * \\
\end{pmatrix}. \]
Thus \eqref{cone_to_take} is of the form
\[ \begin{pmatrix}
-1 & * & \cdots & * & * \\
0 & -1 & \cdots & * & * \\
\vdots & \vdots & \ddots & \vdots & \vdots \\
0 & 0 & \cdots & -1 & * \\
\end{pmatrix} \]
and in particular is split, so the cone on \eqref{cone_to_take} is $FH^{n+1}[1]$ as desired.
\end{proof}
\noindent Thus for example if $H = [-2]$ then $P$ acts on $\im F$ by $[-2n]$, or if $H = [-1]$ then $P$ acts on $\im F$ by $[-n+1]$.

If $F\colon D^b(X) \to D^b(Y)$ is a $\P^n$-functor with $H = [-2]$ then $P$ acts trivially on cohomology, but if $H = [-1]$ or something more interesting then the action on cohomology may be more interesting as well.
\bigskip

Now we can emulate the end of \S\ref{spherical_coho} to show that if $S$ is a K3 surface and $F\colon D^b(S) \to D^b(S^{[n]})$ is a $\P^{n-1}$-functor with $H = [-2]$ then the $\P$-twist $P_F$ associated to $F$ is different from the known autoequivalences of $D^b(S^{[n]})$.  We know that $P_F$ shifts $\im F$ by $-2n+2$ and fixes $\ker R$.  If $\E \in D^b(S)$ is a spherical object with $\E^\perp \ne \emptyset$, $T_\E$ the spherical twist associated to $\E$, and $\varphi\colon \Aut(D^b(S)) \hookrightarrow \Aut(D^b(S^{[n]})$ Ploog's map, then there are non-zero objects which $\varphi(T_\E)$ shifts by $-2n, -2n+1, \dotsc, -2, -1$, and 0, so $P_F$ is not a shift of $\varphi(T_\E)$.  If $\E \in D^b(S^{[n]})$ is a $\P^n$-object with $\E^\perp \ne 0$\footnote{Again, this holds for all known $\E$: if $\E$ is a line bundle, consider $\E \otimes F \O_x^* \otimes F \O_y$, where $x, y \in S$ are distinct points, and if $\E$ is the structure sheaf of a $\P^n \subset S^{[n]}$, consider the structure sheaf of a point not in the $\P^n$.} then $P_\E$ shifts $\E$ by $-2n$ and fixes $\E^\perp$, so $P_F$ is not a shift of $P_\E$.

\subsection{Proof of equivalence}

\begin{thm}
If $F\colon \A \to \B$ is a $\P$-functor then the associated $\P$-twist $P\colon \B \to \B$ is an equivalence.\end{thm}
\begin{proof}
We model our proof on that of Theorem \ref{anno_thm}.  Again we consider the spanning class $\Omega = \ker R \cup \im F$.

First, if $B, B' \in \ker R$ then $\Hom(PB, PB') = \Hom(B, B')$.  Next, if $FA \in \im F$ and $B \in \ker R = \ker L$ then
\begin{align*}
\Hom(PFA, PB)
&= \Hom(FH^{n+1}A[2], B) \\
&= \Hom(H^{n+1}A[2], RB) \\
&= 0 \\
&= \Hom(FA, B) \\[1em]
\Hom(PB, PFA) 
&= \Hom(B, FH^{n+1}A[2]) \\
&= \Hom(LB, H^{n+1}A[2]) \\
&= 0 \\
&= \Hom(B, FA).
\end{align*}
Last, if $FA, FA' \in \im F$ then
\begin{align*}
\Hom(PFA, PFA')
&= \Hom(F H^{n+1} A, F H^{n+1} A') \\
&= \Hom(A, H^{-n-1} RF H^{n+1} A') \\
&= \Hom(A, RF A') \\
&= \Hom(FA, FA').
\end{align*}
Thus $P$ is fully faithful.  To show that $P$ is an equivalence, we show that $\ker P^l = 0$.  Take left adjoints of $PF \cong H^{n+1}F[2]$ to get $H^{-n-1} L[-2] \cong L P^l$, so if $P^l B = 0$ then $H^{-n-1} L B = 0$, so $LB = 0$.  Take left adjoints of the definition of $P$ to get
\[ P^l[1] = \cone(T^l \to F H^{-1} L[-1]), \]
from which we see that if $P^l B = 0$ then $T^l B = 0$, so $B = 0$ as in the proof of Theorem \ref{anno_thm}.
\end{proof}

%% file: cubic_calc.tex

\section{Cubic 4-fold calculation} \label{cubic_calc}

Fix a smooth cubic hypersurface $X \subset \P^5$.  Let $\A \subset D^b(X)$ be Kuznetsov's K3 subcategory
\begin{align*}
\A &:= \langle \O_X, \O_X(1), \O_X(2) \rangle^\perp \\
&= \{ A \in D^b(X) : \RHom(\O_X(i),A)=0\ \mathrm{for}\ i=0,1,2 \},
\end{align*}
let $I\colon \A \to D^b(X)$ be the inclusion, and let $I^l$ and $I^r$ be its left and right adjoints, which exist because $\O_X$, $\O_X(1)$, and $\O_X(2)$ form an exceptional collection.  Kuznetsov \cite{kuz_cubics} has shown that $\A$ satsifies $S_\A = [2]$ and
\[ \HH^i(\A) = \begin{cases}
1 & i = 0 \\
22 & i = 2 \\
1 & i = 4 \\[1ex]
0 & \text{otherwise},
\end{cases} \]
just like the derived category of a K3 surface.

Let $Y \subset \Gr(2,6)$ be the variety of lines on $X$ and let
\[ L = \{ (x,l) \in X \times Y : x \in l \} \]
be the universal line, with projections $X \xleftarrow{q} L \xrightarrow{p} Y$.  We will consider the functor $F := p_* q^* I\colon \A \to D^b(Y)$ and its adjoints $L = I^l q_! p^*$ and $R = I^r q_* p^!$.

In \S\ref{cubic_easy} we explain that $Y$ can be seen as a moduli space of objects in $\A$ and $F$ as the functor induced by the universal object.  In \S\ref{cubic_hard} we prove:
\begin{thm} \label{cubic_thm}
The functor $F\colon \A_X \to D^b(Y)$ is spherical with cotwist $C = [-2]$.  For any point $y \in Y$, the twist $T$ takes $\O_y$ to an object of rank 2, so $T$ is not generated by the previously-known autoequivalences of $D^b(Y)$.  
\end{thm}

\subsection{\texorpdfstring{$Y$}{Y} as a moduli space of objects in \texorpdfstring{$\A$}{A}} \label{cubic_easy}

Let $y \in Y$ be a point, and let $l \subset X$ be the corresponding line.  We wish to describe $R \O_y$.  For general reasons we have
\[ R = S_\A L S_Y^{-1} = L[-2], \]
and thus
\[ R \O_y = L \O_y[-2] = I^l q_! p^* \O_y [-2]. \]
By Grothendieck duality we have $q_!((-)^*) = (q_*(-))^*$, and we check that $(p^* \O_y)^* = p^* \O_y[-4]$ and $p_* q^* \O_y = \O_l$ and $(\O_l)^* = \O_l(1)[-3]$,\footnote{For a refresher on this sort of calculation see \cite[\S3.4]{huybrechts}.} whence
\[ R \O_y = I^l \O_l(1)[-1]. \]
Now $I^l$ is given by left mutation past $\O_X(2)$, $\O_X(1)$, and $\O_X$, where we recall that mutation past $\O_X(i)$ is
\[ \cone\big(\,\O_X(i) \otimes \RHom(\O_X(i),\,-) \longrightarrow \id_X\,\big). \]
We already have $\RHom(\O_X(2), \O_l(1)) = 0$, so mutation past $\O_X(2)$ does nothing.  Mutation past $\O_X(1)$ turns $\O_l(1)[-1]$ into the twisted ideal sheaf $\I_l(1)$.

For mutation past $\O_X$, let $\F_l$ be the ``second syzygy sheaf'' defined by the exact sequence
\[ 0 \to \F_l \to \O_X(-1)^4 \to \O_X \to \O_l \to 0, \]
where for example if $l$ is the line $x_0 = x_1 = x_2 = x_3 = 0$ then the map $\O_X(-1)^4 \to \O_X$ is given by the matrix $\begin{pmatrix} x_0 & x_1 & x_2 & x_3 \end{pmatrix}$.  Then $\F_l$ is a reflexive sheaf of rank 3, locally free away from $l$, and mutation past $\O_X$ turns $\I_l(1)$ into $\F_l(1)[1]$.

Kuznetsov and Markushevich \cite[\S5]{km} have shown that $\F_l$ is a stable sheaf, that $\F_l \not\cong \F_{l'}$ when $l \ne l'$, and that the natural map $T_{y,Y} \to \Ext^1(\F_l, \F_l)$ is an isomorphism.  Thus $Y$ can be seen as a moduli space of objects in $\A$.

In fact if we want $F$ (rather than $R$) to be induced by the ``universal object'' then we should regard $Y$ as the moduli space of objects $\F_l^*(-1)[1]$, which are truly complexes -- they have cohomology sheaves in two degrees -- but this is not a problem.

\subsection{\texorpdfstring{$F$}{F} is spherical} \label{cubic_hard}
To prove Theorem \ref{cubic_thm} we must study the composition
\[ RF = I^r q_* p^! p_* q^* I. \]
The bulk of the work will be in analyzing the middle portion, $q_* p^! p_* q^*$.  Consider the diagram
\[ \xymatrix{
& & L \times_Y L \ar[ld]_{\pi_1} \ar[rd]^{\pi_2} \\
& L \ar[ld]^q \ar[rd]_p & & L \ar[ld]^p \ar[rd]_q \\
X & & Y & & X
} \]
Since $p$ is a $\P^1$-bundle, it is flat, so
\begin{align*}
q_* p^! p_* q^*(-)
&= q_* (\omega_p \otimes p^* p_* q^*(-))[1] \\
&= q_* (\omega_p \otimes \pi_{2*} \pi_1^* q^*(-))[1] \\
&= (q \pi_2)_* (\pi_2^* \omega_p \otimes (q \pi_1)^*(-))[1]
\end{align*}
Thus if we let $\psi = q\pi_1 \times q\pi_2 \colon L \times_Y L \to X \times X$ then the functor $q_* p^! p_* q^*$ is induced by
\[ \psi_* \pi_2^* \omega_p[1] \in D^b(X \times X). \]

\begin{lem} \label{psi_star_whatever}
Let $Z = \im(\psi) \subset X \times X$.  Then
\[ R^i \psi_* \pi_2^* \omega_p = \begin{cases}
I_{\Delta_X/Z}(1,-1) & i = 0 \\
\Delta_* \omega_X & i = 1 \\
0 & \text{otherwise.}
\end{cases} \]
\end{lem}

Before proving this we make two smaller calculations:

\begin{lem}
The variety $Z = \im(\psi) \subset X \times X$ is a complete intersection of two hypersurfaces of bidegrees $(2,1)$ and $(1,2)$ in $X \times X$.
\end{lem}
\begin{proof}
Observe that $Z$ is the closure of the set
\[ \{ (x,y) \in X \times X : \text{$x \ne y$ and the line $\overline{xy}$ lies in $X$} \}. \]
Let $f$ be a polynomial defining $X$, and let $\hat f$ be the polarization of $f$, that is, the unique symmetric trilinear form with $\hat f(v,v,v) = f(v)$ for all $v \in \C^6$.  Given distinct points $x, y \in X$, the line $\overline{xy}$ lies in $X$ if and only if $f(sx+ty) = 0$ for all $s,t \in \C$.  But
\begin{align*}
f(sx+ty) &= s^3\hat f(x,x,x) + 3s^2t \hat f(x,x,y) + 3st^2 \hat f(x,y,y) + t^3 \hat f(y,y,y) \\
&= 3s^2t \hat f(x,x,y) + 3st^2 \hat f(x,y,y)
\end{align*}
since $x,y \in X$, and this vanishes for all $s$ and $t$ if and only if
\begin{align*}
\hat f(x,x,y) &= 0 \\
\hat f(x,y,y) &= 0.
\end{align*}
The first equation has bidegree $(2,1)$; it says that the line is tangent to $X$ at $x$.  The second has bidegree $(1,2)$; it says that the line is tangent at $y$.
\end{proof}

\begin{lem} \label{q_star_O_L}
With $q\colon L \to X$ as above,
\[ R^i q_* \O_L = \begin{cases}
\O_X & i = 0 \\
T_X(-3) & i = 1 \\
0 & \text{otherwise,}
\end{cases} \]
where $T_X$ is the tangent bundle of $X$.\footnote{Thus the general fiber of $q$ is a smooth curve of genus 4.  For the reader's interest, we mention but do not prove that if $X$ is general then $q$ is flat, but for example if $X$ is the Fermat cubic then $q^{-1}(1,-1,0,0,0,0)$ is a smooth cubic surface.

The published version of this paper stated erroneously that $R^1 q_* \O_L = T_X(2)$.  This cannot be correct, because $Rp_* \O_L = \O_Y$ and thus $h^*(\O_L) = h^*(\O_Y) = 1, 0, 1, 0, 1$, whereas $h^*(T_X(2)) = 86, 6, 0, 0, 0$ which would give us $h^*(\O_L) = 1, 86, 6, 0, 0$.  But the correction $R^1 q_* \O_L = T_X(-3) = \Omega^3_X$ resolves this contradiction, because $h^*(\Omega^3_X) = 0, 1, 0, 1, 0$.}
\end{lem}
\begin{proof}
We can naturally embed $L$ in the $\P^3$-bundle $\P T_X$.  Let $\varpi\colon \P T_X \to X$, and let $\O_\varpi(-1)$ denote the tautological rank-1 sub-bundle of $\varpi^* T_X$.  

Let $f$ be a polynomial defining $X$.  Given a point $x \in X$ and a tangent vector $\xi \in T_{x,X}$, we know that $f$ vanishes at $x$, and the derivative $D_\xi f$ vanishes at $x$; moreover the line determined by $x$ and $\xi$ lies in $X$ if and only if the second and third derivatives $D_\xi^2 f$ and $D_\xi^3 f$ vanish at $x$.  Now $D_\xi^2 f$ determines a section of $\O_{\varpi}(2) \otimes \varpi^* \O_X(3)$ which cuts out a hypersurface $M \subset \P T_X$, and $D_\xi^3 f$ determines a section of $(\O_{\varpi}(3) \otimes \varpi^* \O_X(3))|_M$ which cuts out $L$.  Using the exact sequences
\[ 0 \to \O_{\varpi}(-2) \otimes \varpi^* \O_X(-3) \to \O_{\P T_X} \to \O_M \to 0 \]
\[ 0 \to (\O_{\varpi}(-3) \otimes \varpi^* \O_X(-3))|_M \to \O_M \to \O_L \to 0 \]
and the fact that
\[ R^3 \varpi_* \O_{\varpi}(-5) = T_X \otimes \det T_X = T_X(3) \]
we deduce the result.
\end{proof}

\begin{proof}[Proof of Lemma \ref{psi_star_whatever}]
First we argue that
\[ R^0 \psi_* \pi_2^* \omega_p = I_{\Delta_X/Z}(1,-1). \]
Consider the Beilinson resolution of the diagonal $\Delta_L \subset L \times_Y L$:
\[ 0 \to \O_p(-1) \boxtimes_Y \omega_p(1) \to \O_{L \times_Y L} \to \O_{\Delta_L} \to 0, \]
where $\O_p(1)$ is any line bundle on $L$ whose restriction to the fibers of $p$ is $\O_{\P^1}(1)$.  Since $q$ embeds the fibers of $p$ as straight lines in $X$, we can take $\O_p(1) = q^* \O_X(1)$.  Applying $\psi_*$ we get an exact sequence
\[ 0 \to R^0 \psi_* \pi_2^* \omega_p \otimes \O_{X \times X}(-1,1) \to R^0 \psi_* \O_{L \times_Y L} \to \O_{\Delta_X}, \]
where in the first term we have used the projection formula and in the third we have used Lemma \ref{q_star_O_L}.  It is enough now to argue that the middle term is $\O_Z$.  Observe that $\psi$ is an isomorphism away from $\Delta_L$, hence is birational onto its image $Z$.  Next observe that $Z$ is regular in codimension 1, since $Z \setminus \Delta_X \cong (L \times_Y L) \setminus \Delta_L$, and it satisfies Serre's condition S2, being a complete intersection, so it is normal.  Thus Zariski's main theorem gives $R^0 \psi_* \O_{L \times_Y L} = \O_Z$, as desired.
\bigskip

Next we argue that
\[ R^i \psi_* \pi_2^* \omega_p = \begin{cases}
\Delta_* \omega_X & i = 1 \\
0 & i > 1.
\end{cases} \]
Because $\psi$ is an isomorphism away from $\Delta_L$, we see that $R^i \psi_* \pi_2^* \omega_p$ is supported on $\psi(\Delta_L) = \Delta_X$ for $i > 0$, at least set-theoretically, and it is tempting to say that
\begin{equation} \label{tempting}
R^i \psi_* \pi_2^* \omega_p
= R^i \psi_* \left(\pi_2^* \omega_p|_{\Delta_L} \right)
\end{equation}
for $i > 0$.  Assuming that this is true, we have
\begin{align*}
R^i \psi_* \left(\pi_2^* \omega_p|_{\Delta_L} \right)
&= \Delta_* R^i q_* \omega_p \\
&= \Delta_* R^i q_* \omega_L \\
&= \begin{cases}
\Delta_* \omega_X & i = 1 \\ 
0 & i > 1,
\end{cases}
\end{align*}
where in the second line we have used the fact that $\omega_Y = \O_Y$, and in the third we have used Grothendieck duality and Lemma \ref{q_star_O_L}.

To justify \eqref{tempting} we use the theorem on formal functions, which states that the completion of $R^i \psi_* \pi_2^* \omega_p$ along $\Delta_X$ is the inverse limit
\[ \varprojlim_n R^i \psi_* \left(\pi_2^* \omega_p|_{n\Delta_L} \right), \]
where $n\Delta_L$ is the $n^\text{th}$ thickening of $\Delta_L$ in $L \times_Y L$.  We have an exact sequence
\[ 0 \to \O_{\Delta_L}(-n\Delta_L) \to \O_{(n+1) \Delta_L} \to \O_{n \Delta_L} \to 0 \]
for all $n > 0$.  Tensoring with $\pi_2^* \omega_p$ and recalling that $O_{\Delta_L}(-\Delta_L) = \omega_p = \omega_L$, we find that if
\begin{equation} \label{to_kill_with_kodaira}
R^i q_* \omega_L^{n+1} = 0
\end{equation}
for all $i > 0$ and $n > 0$, then the inverse limit stabilizes at the first step, implying \eqref{tempting}.  Since $L$ is the projectivization of the tautological rank-2 sub-bundle on $Y \subset \Gr(2,6)$, we find that $\omega_L = \omega_p = q^* \O_X(-2) \otimes p^* \O_Y(1)$.  Now we can prove \eqref{to_kill_with_kodaira} by direct calculation as in the proof of Lemma \ref{q_star_O_L}, or we can use relative Kodaira vanishing,
since $p^* \O_Y(1)$ is $q$-very ample.
\end{proof}

Now Lemma \ref{psi_star_whatever} gives us an exact triangle
\[ \I_{\Delta_X/Z}(1,-1) \to \psi_* \pi_2^* \omega_p \to \Delta_* \omega_X[-1], \]
which we shift by 1 and rewrite as
\begin{equation} \label{to_compose}
\I_{\Delta_X/Z}(1,-1)[1] \to q_* p^! p_* q^* \to S_X[-4].
\end{equation}
We wish to compose with $I^r$ on the left and $I$ on the right.

\begin{lem} \label{OZ1-1_dies}
$I^r \circ \O_Z(1,-1) \circ I = 0$.
\end{lem}
\begin{proof}
Replace $\O_Z(1,-1)$ with the Koszul complex
\[ \O_{X \times X}(-2,-4) \to \O_{X \times X}(-1,-2) \oplus \O_{X \times X}(0,-3) \to \O_{X \times X}(1,-1). \]
If we apply the functor $I^r \circ \O_Z(1,-1)$ to an object $IA$, where $A \in \A$, then we get a complex
\begin{multline*}
\RGamma(IA(-2)) \otimes I^r(\O_X(-4)) \\
\to \RGamma(IA(-1)) \otimes I^r(\O_X(-2)) \oplus \RGamma(IA) \otimes I^r(\O_X(-3)) \\
\to \RGamma(IA(1)) \otimes I^r(\O_X(-1)).
\end{multline*}
But all the terms vanish: for $i=2,1,0$ we have
\[ \RGamma(IA(-i)) = \RHom(\O_X(i), IA) = 0, \]
and then for the last term we have
\[ I^r(\O_X(-1)) = S_\A I^l S_X^{-1}(\O_X(-1)) = S_\A I^l(\O_X(2))[-4] = S_\A(0) = 0. \qedhere \]
\end{proof}

Now we take the triangle \eqref{to_compose} and apply $I^r$ on the left and $I$ on the right.  By Lemma \ref{OZ1-1_dies}, the first term becomes \[ I^r \circ \O_{\Delta_X} \circ I = I^r I = \id_\A. \]
The second term becomes $RF$.  For the third term, we have 
\[ I^r S_X I [-4] = S_\A [-4] = [-2]. \]
Thus we have an exact triangle
\begin{equation} \label{RF_triangle}
\id_\A \to RF \to [-2].
\end{equation}

Now we can deduce the first statement of Theorem \ref{cubic_thm}, that $F$ is spherical with cotwist $C = [-2]$.  The triangle \eqref{RF_triangle} is split because
\[ \Ext^1([-2], \id_\A) = \Hom(\id_\A, [-3]) = \HH^3(\A) = 0. \]
This argument is valid because we are working with Fourier--Mukai kernels, not just functors; everything is proved rigorously in \cite{kuz_hochschild}.  The first map of \eqref{RF_triangle} agrees with the unit $\eta\colon \id_\A \to RF$ up to a scalar multiple, because
\[ \Hom(\id_\A, RF) = \Hom(\id_A, \id_A \oplus[-2]) = \HH^0(\A) \oplus \HH^{-2}(\A) = \C \oplus 0 \]
and $\eta$ is not zero.  Thus we have $C := \cone \eta = [-2]$, which is an isomorphism, and in \S\ref{cubic_easy} we saw that $R = L[-2] = CL$.
\bigskip

Last we prove the second statement of Theorem \ref{cubic_thm}, that $\rank(T \O_y) = 2$.  We have
\begin{align*}
\rank(T \O_y)
&= \rank(\O_y) - \rank(FR \O_y) \\
&= -\rank(FR \O_y) \\
&= -\chi(FR \O_y, \O_y) \\
&= -\chi(R \O_y, R \O_y),
\end{align*}
where $\chi$ is the Euler pairing.  In \S\ref{cubic_easy} we saw that $R \O_y$ is a shift of a stable sheaf, so $\Ext^{<0}(R \O_y, R \O_y) = 0$ and $\Hom(R \O_y, R \O_y) = \C$.  Since $R \O_y \in \A$ and $S_\A = [2]$, this gives $\Ext^2(R \O_y, R \O_y) = \C$ and $\Ext^{>2}(R \O_y, R \O_y) = 0$.  Finally we have $\Ext^1(R \O_y, R \O_y) \cong T_{y,Y}$, so $\chi(R \O_y, R \O_y) = 1-4+1 = -2$, which proves the claim.

%% file: base_change.tex

\pagebreak
\appendix
\section{Appendix: Cohomology and Base Change}\label{base_change}

The following is well-known to those who know it well, but I could not find a reference.

\newcommand \Dbcoh {D^b_{coh}}
\begin{prop}
Let $X$, $Y$, and $B$ be connected Noetherian schemes of finite type over a field $k$, with $X$ and $Y$ Cohen--Macaulay and $B$ smooth.  Let $f\colon X \to B$ be proper, so that $f_*$ takes $\Dbcoh(X)$ into $\Dbcoh(B)$.\footnote{Recall that our functors are implicitly derived: we mean $Rf_*$, $Lg^*$, etc.}  Let $g\colon Y \to B$ be arbitrary; since $B$ is smooth, $g^*$ takes $\Dbcoh(B)$ into $\Dbcoh(Y)$.  Let $\tilde f$ and $\tilde g$ be as in the diagram
\[ \xymatrix{
X \times_B Y \ar[r]^-{\tilde f} \ar[d]_{\tilde g} & Y \ar[d]^g \\
X \ar[r]_f & B.
} \]
If every irreducible component of $X \times_B Y$ is of the expected dimension $\dim X + \dim Y - \dim B$, then the natural map $g^* f_* \to \tilde f_* \tilde g^*$ is an isomorphism.
\end{prop}
\begin{proof}
Let $\Gamma_f \subset X \times B$ be the graph of $f$ and $\Gamma_g \subset B \times Y$ the graph of $g$.  It is enough to show that on $X \times B \times Y$ we have $\Tor_i(\O_{\Gamma_f \times Y}, \O_{X \times \Gamma_g}) = 0$ for $i > 0$, a condition called ``Tor-independence'' \cite[Thm.~3.10.3]{lip}. Note that $(\Gamma_f \times Y) \cap (X \times \Gamma_g) \cong X \times_B Y$.

First I claim that $\Gamma_f$ is locally cut out of $X \times B$ by a regular sequence.  Since $B$ is smooth, the diagonal $\Delta \subset B \times B$ is locally cut out by a regular sequence of $n$ functions, where $n = \dim B$.  
Thus $\Gamma_f = (f \times 1)^{-1}\Delta$ is locally cut out by $n$ functions, which \emph{a priori} may not be a regular sequence; but $X \times B$ is Cohen--Macaulay \cite{ty}, so a sequence of $n$ functions is regular if and only if it cuts out a subscheme of codimension $n$ \cite[Thm.\ 17.4(iii)]{mat}, and the codimension of $\Gamma_f \cong X$ is indeed $n$.

Thus $\Gamma_f \times Y$ is locally cut out of $X \times B \times Y$ by a regular sequence of $n$ functions, so locally we can resolve $\O_{\Gamma_f \times Y}$ by a Koszul complex.  Tensoring with $\O_{X \times \Gamma_g}$, we see that the higher Tors vanish if the sequence remains regular when restricted to $X \times \Gamma_g$.  Since $X \times \Gamma_g \cong X \times Y$ is Cohen--Macaulay and the subscheme $(\Gamma_f \times Y) \cap (X \times \Gamma_g) \cong X \times_B Y$ cut out by the restricted sequence has codimension $n$ by hypothesis, we are done.
\end{proof}

The dimension hypothesis is necessary:  Let $B$ be a smooth surface, $Y$ a point, and $X$ the blowup of $B$ at $g(Y)$.  Then $E := X \times_B Y$ is the exceptional line, whose dimension is $1 > 0 + 2 - 2$, and we find that $g^* f_* \tilde g_* \O_E(-1) = 0$ while $\tilde f_* \tilde g^* \tilde g_* \O_E(-1) = \O_Y[1]$.

The smoothness of $B$ is necessary:  Let $B$ be the cone $xy=z^2$ in $\mathbb A^3$, let $X$ be the line $y=z=0$, let $Y$ be the line $x=z=0$, and let $f$ and $g$ be the inclusions, so $X \times_B Y = X \cap Y$ is the origin, which is of the expected dimension.  Using the resolution
\[ \dotsb \to \O_B^2
\xrightarrow{\left(\begin{smallmatrix} y & -z \\ -z & x \end{smallmatrix}\right)} \O_B^2
\xrightarrow{\left(\begin{smallmatrix} x & z \\ z & y \end{smallmatrix}\right)} \O_B^2
\xrightarrow{\left(\begin{smallmatrix} y \\ -z \end{smallmatrix}\right)} \O_B
\to f_*\O_X \to 0 \]
we find that $\Tor_i(f_*\O_X, g_*\O_Y) = \O_\text{origin}$ for all $i \ge 0$.  Thus $g_* g^* f_* \O_X = f_*\O_X \otimes g_*\O_Y$ is different from $g_* \tilde f_* \tilde g^* \O_X = g_* \tilde f_* \O_{X \cap Y}$.

The Cohen--Macaulay hypothesis is also necessary, as we see from the following example based on \cite{speyer_mo}.\footnote{There are simpler examples, but this is the simplest one I know in which all the spaces are normal.}  Let $B =  \mathbb A^6$.  Let $C$ be the Fermat cubic curve $\{ x_0^3 + x_1^3 + x_2^3 = 0 \} \subset \P^2$, and let $X \subset \mathbb A^6$ be the affine cone over $C \times \P^1 \subset \P^2 \times \P^1 \subset \P^5$; this is not Cohen--Macaulay since the surface $C \times \P^1$ is not arithmetically Cohen--Macaulay (it has $H^1(\O) = 1$).  Let $Y$ be a generic 2-plane through the origin in $\mathbb A^6$, and let $g$ be the inclusion.  With Macaulay2 \cite{M2} we calculate that $X \cap Y$ is a scheme of length 7 supported at the origin, which is of the expected dimension, but $\Tor_1(f_*\O_X, g_*\O_Y) = \O_\text{origin}$:
\begin{verbatim}
R = QQ[x_0..x_2, y_0,y_1] -- P^2 x P^1
S = QQ[z_0..z_5] -- P^5 or A^6
segre = map(R, S,
  {x_0*y_0, x_1*y_0, x_2*y_0, x_0*y_1, x_1*y_1, x_2*y_1})
IX = preimage_segre ideal(x_0^3 + x_1^3 + x_2^3)
IY = ideal random(S^{1}, S^3) -- three random linear forms
dim(IX + IY) -- answer is 0
degree(IX + IY) -- answer is 7
length Tor_1(comodule IX, comodule IY) -- answer is 1
\end{verbatim}
So again $g_* g^* f_* \O_X \ne g_* \tilde f_* g^* \O_X$.